\documentclass[10pt]{amsart}

\usepackage{hyperref}
\usepackage{graphicx}
\usepackage{tikz}
\usepackage{amsthm}
\usepackage{multicol}

\usepackage{amssymb}
\usepackage{amsmath}

\usepackage{ulem}

\let\newpf\proof \let\proof\relax

\def\bm{\begin{matrix}}
\def\em{\end{matrix}}

\newcommand{\bt}{\begin{thm}}
\newcommand{\et}{\end{thm}}

\newcommand{\bl}{\begin{lemma}}
\newcommand{\el}{\end{lemma}}

\newcommand{\beq}{\begin{eqnarray}}
\newcommand{\eeq}{\end{eqnarray}}

\def\be{\begin{equation}}
\def\ee{\end{equation}}

\def\ba{{\begin{align}}}
\def\ea{{\end{align}}}

\def\0{{\mathbf 0}}

\def\cal{\mathcal}

\newtheorem{thm}{Theorem}[section]
\newtheorem{cor}[thm]{Corollary}
\newtheorem{lem}[thm]{Lemma}
\newtheorem{lemma}[thm]{Lemma}
\newtheorem{claim}[thm]{Claim}
\newtheorem{prop}[thm]{Proposition}

\theoremstyle{remark}
\newtheorem{rem}{Remark}[section]

\numberwithin{equation}{section}

\def \bn {\hfill \\ \smallskip\noindent}

\theoremstyle{definition}
\newtheorem{defn}{Definition}[section]

\def\proof{\bn {\bf Proof.} }

\def\note#1
{\marginpar

{\tiny $\leftarrow$
\par
\hfuzz=20pt \hbadness=9000 \hyphenpenalty=-100 \exhyphenpenalty=-100
\pretolerance=-1 \tolerance=9999 \doublehyphendemerits=-100000
\finalhyphendemerits=-100000 \baselineskip=6pt
#1}\hfuzz=1pt}

\newcommand{\C}{{\mathbb C}}

\newcommand{\E}{{\mathbb E}}

\newcommand{\N}{{\mathbb N}}
\newcommand{\Q}{{\mathbb Q}}
\newcommand{\R}{{\mathbb R}}
\newcommand{\T}{{\mathbb T}}

\newcommand{\Z}{{\mathbb Z}}

\def\B0{{\bold{0}}}

\catcode`\@=12

\def\Empty{}
\newcommand\oplabel[1]{
  \def\OpArg{#1} \ifx \OpArg\Empty {} \else
  	\label{#1}
  \fi}

\newcommand{\comm}[1]{}
\newcommand{\comment}[1]{}

\begin{document}

\title[Spectral Dimension for singular Jacobi operators]{Spectral Dimension for $\beta$-almost periodic singular Jacobi operators and the extended Harper's model }

\author{Rui Han, Fan Yang and Shiwen Zhang}


\thanks{}

\begin{abstract}
We study fractal dimension properties of  singular Jacobi operators. We
prove quantitative lower spectral/quantum dynamical bounds for general 
operators with strong repetition properties and controlled 
singularities. 
For analytic quasiperiodic
Jacobi operators in the positive Lyapunov exponent regime, we obtain a 
sharp arithmetic
criterion of full spectral dimensionality. The applications include the
extended Harper's model where we obtain arithmetic results on spectral 
dimensions and quantum dynamical exponents.
\end{abstract}

\maketitle

\section{Introduction}
In this paper, we study self-adjoint Jacobi operators on ${\ell}^2(\Z)$ given by:
\begin{equation}\label{def:Jacobi}
(Hu)_n=w_n u_{n+1} +\overline{w}_{n-1} u_{n-1}+v_n u_n,\ n\in\Z
\end{equation}
where $w_n \in \mathbb{C} \setminus \{ 0 \}$ and $v_n \in \R$ are bounded  sequences in $n$. 
If $w_n\equiv 1$, $H$ is a discrete Schr\"odinger operator. 
We will focus on singular Jacobi operators, where the off-diagonal sequence $w_n$ has an accumulation point $0$ at $\pm\infty$. 
A prime example of such operator, both in math and in physics literature,
is the extended Harper's model (EHM), see (\ref{def:EHM}). 


We are interested in the fractal decomposition of the spectral measure and quantitative spectral/quantum dynamical bounds.
In a recent work of Jitomirskaya and Zhang \cite{JZ}, many quantitative criteria of fractal dimensions of spectral measure for lattice Schr\"odinger operators were obtained. 
Their criteria have various applications to quasiperiodic Schr\"ondiger operators, e.g., the almost Mathieu operator or the Sturmian Hamiltonians. 
However, whether their results could be applied to the Jacobi case, in particular the singular Jacobi case, has remained a question.
Indeed, many generalizations of the spectral theory from the Schr\"odinger case to the singular Jacobi case have shown to be highly nontrivial, e.g. \cite{JKS,JMarx12CMP,M14,HJ,H18,AJM,Y}.






In this paper we give general sufficient conditions for spectral continuity in the singular Jacobi case, see Theorem \ref{thm:speconti}.
We show spectral continuity follows if the parameters $w_n,v_n$ of $H$ satisfy: (i) strong repetition properties; (ii) control of the averaged closeness between $w_n$ and 0. 
In particular, condition (ii) is imposed to control the strength of singularity.
We also show such operators exist widely in the general context of quasiperiodic setting. 
In the positive Lyapunov exponent regime of analytic quasiperiodic Jacobi operators, the general statement leads to the first arithmetic if-and-only-if criterion for full spectral dimensionality. 
Notably, our results have applications to the extended Harper's model 
 for both spectral and quantum dynamical properties.

Our proof is based on a general dynamical system approach, which has recently shown to be extremely powerful in the study of spectral properties, e.g. \cite{L1,BGT,KKL,GKT,DT,JL1,LS99}.
The eigenvalue equation of \eqref{def:Jacobi}, is associated to a linear cocycle system, see \eqref{eq:eigeneq}. 
In the Schr\"odinger case the cocycles are $\mathrm{SL}(2,\R)$-valued, whereas in the singular Jacobi case the cocycles are $\mathrm{GL}(2,\C)$-valued with determinants approaching zero along a sub-sequence.
This presents the main obstruction in \cite{JKS,JMarx12CMP,M14,HJ,H18,AJM,Y} and in our paper.

It was shown in \cite{JZ} that the fractal dimensions of spectral measures depend on the competition between the quality of repetitions and the growth of the Schr\"odinger cocycles.
Such competition was resolved in the $\mathrm{SL}(2,\R)$ setting involving delicate algebraic arguments, which are difficult to carry out directly in the $\mathrm{GL}(2,\C)$ setting due to the presence of singularity. 
To reduce 
to $\mathrm{SL}(2,\R)$ case, we employ a family of conjugacies which were first introduced in a recent work of Avila-Jitomirskaya-Marx \cite{AJM}.
Such regularization moves the singularity into the conjugate matrices.
The main technical accomplishment of our work is to develop general quantitative estimates (see Lemmas \ref{lem:rT}, \ref{lem:normA-A} and \ref{lem:Abeta}) of the conjugacy under assumptions (i) and (ii).
The successful combination of these estimates with the mechanism in \cite{JZ} proves the quantitative spectral continuity results for the singular Jacobi case. 

We show the assumptions (i) and (ii) hold for singular Jacobi operators over a quasiperiodic base.
In particular, the proof of (ii) is close in spirit to the characterization of singularity in \cite{HJ} (see also \cite{JLiuMM, JY}). 
Here we need to study the finer decomposition of the singular spectral measure, thus a strengthened characterization is developed.
Moreover, our estimates hold for general $C^k$ sampling functions with finitely many non-degenerate zeros, which reduces the analytic regularity requirements in \cite{JLiuMM,JY,HJ}.
This part is also of independent interest in the study of uniform upper-semi continuity of the Lyapunov growth.

The rest of this paper is organized in the following way. In section \ref{sec:result}, we give all the definitions and state our main results. After giving the preliminaries
in section \ref{sec:prelim}, we proceed to discuss the $(\Lambda,\beta)$ bound in the quasiperiodic case in section \ref{sec:LambdaBound}. In section \ref{sec:spe-conti}, we prove the general spectral continuity results. In section
\ref{sec:speSingular}, we focus on the analytic quasiperiodic Jacobi operator and prove arithmetic if-and-only-if criterion 
for full spectral dimensionality. In the last section, we discuss the explicit parameter partitions for the extended Harper's model.

\section{Main results}\label{sec:result}


To formulate the main results, we introduce the following definitions. 
\begin{defn}\label{def:beta-almost}
	A sequence $\{a_n\}_{n\in\Z}$ is said be to $\beta$-$q$ almost periodic if
	there exist $\delta >0$, $\beta>0$, $q\in \N$, such that the following holds:
	\begin{align}
	\max_{|m|\leq e^{\delta\beta q}} |a_{m} -a_{m\pm q}|
	\leq {e^{-\beta q}}. \label{eq:a-beta}
	\end{align}
	We say $\{a_n \}_{n\in \Z}$ is $\beta$-almost periodic (about $q_n$) if
	there exists a sequence of positive integers $q_n\to\infty$, such that $\{a_n\}$ is $\beta$-$q_n$ almost periodic.
\end{defn}

\begin{rem}
	The $\beta$-almost periodicity was first introduced in \cite{JZ} to study quantitative spectral bounds in the Schr\"odinger case. 	
	Note that the $\beta$-almost periodicity does not imply the almost periodicity in the usual sense. A typical example is the sequence generated by skew-shift map $(x,y)\mapsto (x+y,y+2\alpha)$ with a smooth sampling function $f(x,y)$ on $\T^2$. The sequence $v_n=f(x+ny+n(n-1)\alpha, y+2n\alpha)$ is $\beta$-almost periodic for typical $\alpha$, but not almost periodic for any $\alpha$.
	\end{rem}

\begin{defn}\label{def:Lambda-bound}
	We say $w_n$ is $(\Lambda,\beta)$-$q$ bounded if there exist $\Lambda>0, \beta>0,\delta>0$, and $q\in \N$,
	such that
	\begin{align}\label{eq:Lambda-bound}
	\min_{|m|\le e^{\delta\beta\,q}}\prod_{j=m}^{m+q-1}   |w_j|>e^{-\Lambda\,q}.
	\end{align}
	We say $w_n$ is $(\Lambda,\beta)$ bounded (about $q_n$) if there exists a sequence of positive integers $q_n\to\infty$, such that $w_n$ is $(\Lambda,\beta)$-$q_n$ bounded.
\end{defn}
\begin{rem}\label{rmk:Gordon}	If we only consider the maximum of (\ref{eq:a-beta}) and \eqref{def:Lambda-bound} over $|m|\le 2q_n$, then the standard Gordon-type argument will be enough to show the absence of point spectrum for the associated Jacobi operator, provided $\beta\gtrsim\Lambda$. Assume further the Lyapunov exponent  is positve, then the operator has purely singular continuous spectrum by Kotani theory \cite{K}, see e.g. \cite{D07,BHJ}.   See more discussion on the Gordon-type argument and purely singular continuous spectrum in \cite{D17} and references therein. \end{rem}

Let  $\mu$  be the spectral measure of the Jacobi operator given as in (\ref{def:Jacobi}). The fractal properties of $\mu$ are closely related to the  
boundary behavior of its Borel transforms, see e.g. \cite{DJLS}. Let
\begin{equation}\label{def:Mfunction}
    M(E+ i\varepsilon)=\int\frac{{\rm d}\mu(E')}{E'-(E+ i\varepsilon)}
\end{equation}
be  the (whole line) Weyl-Titchmarsh $m$-function of $H$. We are interested in the following fractal dimension of $\mu$:
\begin{defn}\label{def:dimspe} 
We say $\mu$ is (upper) $\gamma$-spectral continuous if for some $\gamma\in(0,1)$ and $\mu$ a.e. $E$, we have
\begin{equation}\label{def:speconti}
    \liminf_{\varepsilon\downarrow0}\varepsilon^{1-\gamma}|M(E+ i\varepsilon)|<\infty.
\end{equation}
Define the (upper) spectral dimension  of $\mu$ to be
  \begin{align}\label{eq:dimspe}
 {\rm dim}_{\rm spe}(\mu)=\sup\big\{\gamma\in(0,1):\ \mu\ \textrm{is}\ \gamma\textrm{-spectral continuous}\big\}.
 \end{align}
 
\end{defn}

Our first result is about  spectral continuity and the lower bound on the spectral dimension. 
\begin{thm}\label{thm:speconti}
Let $H$ be given as in (\ref{def:Jacobi}) and let $\mu$ be the spectral measure of $H$.  Assume that there are positive constants $\Lambda,\beta,\delta$ and a sequence of positive integers $q_n\to\infty$ such that $w_n,v_n$  are $\beta$-almost periodic and $w_n$ is  $(\Lambda,\beta)$ bounded about $q_n$. There exists an explicit constant $C=C(\delta,\Lambda,\|w\|_{\infty},\|v\|_{\infty})>0$, such that if $\beta>C$ and $\gamma< 1-\frac{C}{\beta}$, then 
$\mu$ is $\gamma$-spectral continuous. Consequently, we have the following lower bound on the spectral dimension of $\mu$: 
\begin{align}\label{eq:spe-lower}
{\rm dim}_{\rm spe}(\mu) \geq 1-\frac{C}{\beta}.
\end{align}
\end{thm}

%

We will formulate a more precise lower bound (specifying the dependence
of $C$ on $\delta,\Lambda,\|w\|_{\infty},\|v\|_{\infty}$) in Theorem \ref{thm:trace}.


It is well known that periodicity implies absolute continuity. We actually prove a quantitative weakening version of this result: $\beta$-almost periodicity implies $\gamma$-spectral continuity. On the other hand, it is well known that Gordon condition implies absense of point spectrum, which predicts purely singular continuous spectrum in many situations. Our result distinguishes the singular continuous spectrum further according to their spectral dimensions. This can be viewed as a quantitative
strengthening of Gordon-type results. 
Quantitative results directly linking easily formulated properties of the potential to dimensional/quantum dynamical results were first proved in \cite{JZ} for the Schr\"odinger case. 
Theorem \ref{thm:speconti} was a further generalization of this type of estimates to more general singular Jacobi operators. \\


An important context where we have generic $\beta$-almost periodicity and $(\Lambda,\beta)$ bound is the quasiperiodic Jacobi operators with smooth sampling functions defined as follows. Consider real and complex valued sampling functions $v:\T\mapsto \R$ and  $c:\T\mapsto \mathbb{C}$. We also assume $\ln|c|\in L^1(\T)$, which is the minimum requirement for the Lyapunov exponent to exist.  Let $H_{\alpha,\theta}=H_{\alpha,\theta,c,v}$ be the Jacobi operator on ${\ell}^2(\Z)$ given by:
\begin{align}\label{def:qpJacobi}
(H_{\alpha,\theta}u)_n =c(\theta+n\alpha) u_{n+1} + \bar{c}\big(\theta+(n-1)\alpha\big)u_{n-1} +v(\theta+n\alpha) u_n,\ n\in\Z,
\end{align}
where  $\theta\in\T:=[0,1]$ is the phase, $\alpha\in [0,1]\backslash\Q$ and $\bar c(\theta)$ is the complex conjugate of $c(\theta)$ in the usual sense.  

Given $\alpha$, let ${p_n}/{q_n}$ be the continued fraction approximants to $\alpha$. Define
\begin{equation}\label{def:beta}
    \beta(\alpha):=\limsup_n\frac{\ln q_{n+1}}{q_n}\in[0,\infty].
\end{equation}

It is easy to check that for any Lipschitz continuous sampling functions $v$ and $c$, the sequences $v(\theta+n\alpha),c(\theta+n\alpha)$ are $\beta$-almost periodic as defined in (\ref{eq:a-beta}) for any $\theta\in\T$ and any $\beta<\beta(\alpha)/2$. See a proof about this simple fact in section \ref{sec:LambdaBound}. Furthermore, if we require some non-degenerate regularity near the zeros of $c$,  then $c(\theta+n\alpha)$ will be $(\Lambda,\beta)$ bounded for a.e. $\theta$. As a consequence of Theorem \ref{thm:speconti}, we have spectral continuity for a.e. $\theta$ for (\ref{def:qpJacobi}). More precisely,  
let $\mu_{\alpha,\theta}$ be the
spectral measure of $H_{\alpha,\theta}$ (\ref{def:qpJacobi}), we have: 

\begin{cor}\label{cor:LipJacobi} Assume $v(\theta)$ is Lipschitz continuous on $\T$ and $c(\theta)$ is $C^{k}$ continuous on $\T$ with finitely many non-degenerate zeros
	\footnote{We say $\theta_0\in\T$ is a non-degenerate zero of $f\in C^{k}(\T,\C)$  if $f(\theta_0)=0$ and $f^{(k)}(\theta_0)\neq0 $.}. For all $k=1,2,\cdots$, there exists an explicit constant $C=C(c,v,k)>0$ and a full measure set $\Theta=\Theta(\alpha,c)\subsetneq\T$, only depending on $\alpha$ and the zeros of $c(\theta)$ with the following properties: suppose  $\beta(\alpha)>C$, then for any $\theta\in\Theta$, 
\begin{description}
  \item[(a)] $H_{\alpha,\theta}$ has no eigenvalues in the spectrum;
  \item[(b)] the spectral dimension of $\mu_{\alpha,\theta}$ is bounded from below as:
  \begin{align}\label{eq:spe-lower-Lip}
     {\rm dim}_{\rm spe}(\mu_{\alpha,\theta}) \geq 1-\frac{C}{\beta}.
  \end{align}
  In particular, if $\beta(\alpha)=\infty$, then for a.e. $\theta$, ${\rm dim}_{\rm spe}(\mu_{\alpha,\theta})=1$. 
\end{description}
\end{cor}

We will prove $(\Lambda,\beta)$ bound of $c(\theta+n\alpha)$ for a.e. $\theta$ in section \ref{sec:LambdaBound} and then part (b) follows directly from Theorem \ref{thm:speconti}. The main ingredient is one fundamental estimate (see Lemma \ref{lem:AJ09}) about the trigonometric product over irrational rotation in \cite{AJ09}.  Similar arguments have been used in \cite{JLiuMM,HJ,JY} to study the arithmetic criterion of purely singular continuous spectrum.  In those papers, the authors considered periodic approximation based on Gordon-type arugment. The growth of the transfer matrix only need to be controlled  within at most two periods. In our case, the quantitative spectral continuity relies on $(\Lambda,\beta)$ bound over exponentially many periods. The use of Lemma \ref{lem:AJ09} is more delicate and involved. See more details in Lemma \ref{lem:LipBound}. 

As mentioned before, the absence of point spectrum in part (a) is a direct consequence of the $(\Lambda,\beta)$ boundedness of $c(\theta+n\alpha)$ and the standard Gordon-type argument. In view of Definition \ref{def:dimspe}, it is easy to check that point measure has zero (spectral/Haursdoff/packing) dimension. (\ref{eq:spe-lower-Lip}) implies that the spectral measure $\mu_{\alpha,\theta}$ has positive spectral dimension for $\beta>C$. Part (a) can also be derived as a corollary of (\ref{eq:spe-lower-Lip}).  An interesting question that remained here is whether the assumption on $c(\theta)$ can be weakened: For example, could any Lipschitz continuous function with finitely many zeros generate a  $(\Lambda,\beta)$ bounded sequence?  Will the associated Jacobi operator have absence of point spectrum and full spectral dimension? 
We will not go further in this direction in the current paper. We are planning to answer some of these questions in another paper (under preparation). \\

It is clear that our general results (\ref{eq:spe-lower}) and (\ref{eq:spe-lower-Lip}) only go in one direction, as
even absolute continuity of the spectral measures does not imply
$\beta$-almost periodicity for $\beta>0.$
However, in the important context of analytic quasiperiodic  operators (e.g. EHM) this leads to a sharp
if-and-only-if result in the positive Lyapunov exponent  regime. 

Let $H_{\alpha,\theta}$ be the Jacobi operator on ${\ell}^2(\Z)$ defined as in (\ref{def:qpJacobi}). The Lyapunov exponent of  $H_{\alpha,\theta}$ at energy $E$ is defined through the associated skew-product over irrational rotations (quasiperiodic cocycles). For any irrational $\alpha$, the Lyapunov exponent is only a function of $E,\alpha$ and is independent of $\theta$, therefore, denoted as $L(E,\alpha)$. See more basic properties and discussions about Lyapunov exponent in section \ref{sec:prelim}.

 Assume further $v,c$ of $H_{\alpha,\theta}$ are analytic on $\T$ with real and complex values, respectively.  Let  $\mu_{\alpha,\theta,\Sigma_+}$ be the restriction of the spectral measure $\mu_{\alpha,\theta}$ of $H_{\alpha,\theta}$ on $\Sigma_+:=\{E\in \sigma(H_{\alpha,\theta}): L(E,\alpha)>0\}$.  We have the following sharp estimate on the spectral dimension of $\mu_{\theta,\alpha,\Sigma_+}$:

\begin{thm}\label{thm:dimspe=1}  For any $\alpha\in[0,1]$, let $\beta(\alpha)$ be defined as in (\ref{def:beta}). For any analytic sampling functions $v$ and $c$, 
there is a full Lebesgue measure set $\Theta=\Theta(\alpha,c)\subset\T$ explicitly depends on $\alpha$ and $c$ such that  for any $\theta\in\Theta$,  ${\rm dim}_{\rm spe}\big(\mu_{\alpha,\theta,\Sigma_+}\big)=1$  if and
only if $\beta(\alpha)=\infty.$
\end{thm}

The proof of Theorem \ref{thm:dimspe=1} contains two parts. Clearly, the `if' part of Theorem \ref{thm:dimspe=1} is a direct consequence of spectral continuity and follows from Corollary \ref{cor:LipJacobi}. The `only if' part is usually referred to as the so-called spectral singularity, defined through the singular boundary behavior of the $m$ function. More precisely, we say the spectral measure $\mu$ is (upper) $\gamma$-spectral singular  if for some  $\gamma\in(0,1)$ and $\mu$ a.e. $E$,
\begin{equation}\label{def:speSingular}
    \liminf_{\varepsilon\downarrow0}\varepsilon^{1-\gamma}|M(E+ i\varepsilon)|=+\infty. 
\end{equation}
Define
\begin{equation}\label{def:spedim-singular}
    {\rm dim}_{\widetilde {\rm spe}}=\inf\big\{\gamma\in(0,1):\ \mu\ \textrm{is}\ \gamma\textrm{-spectral singular}\big\}.
\end{equation}
Obviously, ${\rm dim}_{ {\rm spe}}\le {\rm dim}_{\widetilde {\rm spe}}$. Theorem \ref{thm:dimspe=1} also holds for ${\rm dim}_{\widetilde {\rm spe}}$.  We actually can prove the following local quantitative upper bound of the spectral dimension which completes the sufficient part of Theorem \ref{thm:dimspe=1}. 

\begin{thm}\label{thm:speSingular}
Consider the quasiperiodic Jacobi operator defined in (\ref{def:qpJacobi}) with analytic sampling functions $v,c$. Let $L(E)$ be the associated Lyapunov exponent defined in (\ref{def:Lyp}). Assume $L(E) \geq a >0$ on a compact set $S$. Consider the spectral measure $\mu_{\alpha,\theta}$ restricted on $S$, denoted by $\mu_{\alpha,\theta,S}$. Suppose $\beta(\alpha)<\infty$, then there is a $C=C(a, v, c, S)>0$ and a full Lebesgue measure set $\Theta=\Theta(\alpha,c)$, such that for any $\theta\in\Theta$, $\mu_{\alpha,\theta}$ is $\gamma$-spectral singular  for any $\gamma\ge \frac{1}{1+\frac{C}{\beta}}$. Consequently, 
\begin{equation}\label{eq:spe-upper}
\rm{dim}_{spe} (\mu_{\alpha,\theta,S})\le \rm{dim}_{\widetilde {\rm spe}} (\mu_{\alpha,\theta,S}) \leq \frac{1}{1+\frac{C}{\beta}}<1. 
\end{equation}
\end{thm}

The spectral singularity can be viewed as ``weak-type of localization". It involes the decay of the Green's function in a finite box with a low density (see Lemma \ref{lem:density}). Such decay/localization density was previously known either with a strong non-resonance condition on $\omega$ (e.g. $\beta(\omega)=0$, see \cite{DT}), or for a concrete example with $\beta(\omega)\lesssim L$ (see \cite{AYZ,JLAMO}). Such a phenomenon was first found in \cite{JZ} for general analytic quasiperiodic Schr\"odinger operators with extremetely large $\beta$. Two crucial ingredients for the quantitative spectral singularity are:
\begin{enumerate}
	\item quantitative subordinate theory (Jitomirskaya-Last inequality, Lemma \ref{lem:JLineq});
	\item existence of generalized eigenfunctions with sub-linear growth by Last-Simon estimate (Lemma \ref{lem:ls99}). 
\end{enumerate}
 Theorem \ref{thm:speSingular} generalizes the result for Schr\"odinger operators in \cite{JZ} to singular Jacobi operators. The techniques to deal with the singular Jacobi case are more involved and very delicate in view of the quantitative estimates (\ref{eq:spe-upper}). In section \ref{sec:speSingular}, we reduce the proof of 
Theorem \ref{thm:speSingular} to a quantitative result (see Lemma \ref{lem:JZ3.1}) obtained in \cite{JZ}. 
One key observation in \cite{JZ} is that the norm of the analytic transfer matrix can be approximated by trigonometric polynomials with uniform linear degree. The generalization of this result to the meromorphic transfer matrix in our case (see Lemma \ref{lem:Deltan}) becomes an important part of the proof of Theorem \ref{thm:speSingular}. 





\subsection{Applications to the extended Harper's model.}\label{sec:EHM}
Quasiperiodic Jacobi operators arise naturally from the study of tight-binding electrons on a two-dimensional lattice exposed to a perpendicular magnetic field. A more general model is the extended Harper's model (EHM), defined as follow:
\begin{equation}\label{def:EHM}
(H_{\lambda,\alpha,\theta}u)_n=c_{\lambda}(\theta+n\alpha)u_{n+1}+\bar{c}_{\lambda}\big(\theta+(n-1)\alpha\big)u_{n-1}+2\cos{2\pi (\theta+n\alpha)} u_n,\ n\in\Z. 
\end{equation}
Here,  
\begin{align} \label{def:clambda}
c_{\lambda}(\theta)=\lambda_1 e^{-2\pi i(\theta+\frac{\alpha}{2})}+\lambda_2 +\lambda_3 e^{2\pi i(\theta+\frac{\alpha}{2})},
\end{align}
$\bar c_{\lambda}(\theta)$ is the complex conjugate of $c_{\lambda}(\theta)$ in the usual sense and  $\lambda=(\lambda_1,\lambda_2,\lambda_3)\in\R^3$ are real coupling constants. EHM was introduced by D.J.Thouless in 1983 \cite{T83}, which includes the AMO as a special case.

The extended Harper's model is a prime example of quasiperiodic Jacobi
matrix. It has attracted great attention from both mathematics and physics (see e.g. [7, 10, 20])
literature in the past several decades. 
Recent developments on the spectral theory of the AMO and EHM include:
pure point spectrum for Diophantine frequencies in the positive Lyapunov exponent regime $\mathrm{I}^{\mathrm{o}}$ \cite{JKS};
explicit formula for the Lyapunov exponent $L(E,\lambda)$ (see (\ref{def:EHMLyp})) on the spectrum throughout all the three regions \cite{JMarx12CMP};
dry ten Martini problem for Diophantine frequencies in the self-dual regions \cite{H18};
complete spectral decomposition for all $\alpha$ and a.e. $\theta$ in the zero Lyapunov exponent regiems \cite{AJM};
and arithmetic spectral transition in $\alpha$ in the positive Lyapunov exponent regime \cite{HJ}.

As a central example of the analytic quasiperiodic singular Jacobi operators, Theorem \ref{thm:dimspe=1} can be applied to the extended Harper's model $H_{\lambda,\alpha,\theta}$ defined in (\ref{def:EHM}). As a consequence of the Lyapunov exponent formula of EHM in terms of the coupling constants $\lambda=(\lambda_1,\lambda_2,\lambda_3)$, we have more explicit conclusions on the full spectral dimensionality of EHM.  Moreover,  
our lower bounds in Theorem \ref{thm:dimspe=1} are effective for $\beta>\max\{C\sup_{E\in\sigma(H)}L(E),0\}$ by some simple scaling argument (see Lemma \ref{lem:LipBound} and section \ref{sec:EHM}). Thus the
range of $\beta$  is increased  for smaller Lyapunov exponents. In
particular, we obtain full spectral dimensionality as long as $\beta(\alpha)>0,$ when Lyapunov
exponents are zero on the spectrum. This applies, in particular, to the critical EHM.

Consider the following three parameter regions of $\lambda=(\lambda_1,\lambda_2,\lambda_3)\in\R^3$:
\begin{align*}
&{\cal R}_1=\left\{\lambda\in\R^3:\, 0<\lambda_1+\lambda_3<1,0<\lambda_2<1 \right\}.\\
&{\cal R}_2=\left\{\lambda\in\R^3:\, \lambda_2>\max\{\lambda_1+\lambda_3,1\},\lambda_1+\lambda_3\ge0 \ \ 
{\rm or}\ \ \lambda_1+\lambda_3>\max\{\lambda_2,1\},\lambda_1\neq\lambda_3,\lambda_2>0 \right\} .\\
&{\cal R}_3=\{\lambda\in\R^3:\, 0\le \lambda_1+\lambda_3\le 1,\lambda_2=1 \ \  {\rm or}\ \ \lambda_1+\lambda_3\ge \max\{\lambda_2,1\},\lambda_1=\lambda_3,\lambda_2>0 \} .
\end{align*}

\begin{cor}\label{cor:EHM}
	Let $\mu_{\lambda,\alpha,\theta}$ be the spectral measure of EHM: $H_{\lambda,\alpha,\theta}$. For any $\alpha\in[0,1]$, there is a full measure set $\Theta=\Theta(\alpha)\subset\T$ such that for all $\theta\in\Theta$, the following hold:
	\begin{enumerate}
		\item For $\lambda\in {\cal R}_1$,  ${\rm dim}_{\rm spe}\big(\mu_{\lambda,\alpha,\theta}\big)=1$  if and
		only if $\beta(\alpha)=\infty$.
		
		\item For $\lambda\in {\cal R}_2$ and for all $\alpha\in[0,1]$, ${\rm dim}_{\rm spe}\big(\mu_{\lambda,\alpha,\theta}\big)=1$.
		
		\item For $\lambda\in {\cal R}_3$, ${\rm dim}_{\rm spe}\big(\mu_{\lambda,\alpha,\theta}\big)=1$  if $\beta(\alpha)>0$.
		
	\end{enumerate}
\end{cor} 
We will see the explicit formula of the Lyapunov exponent and the spectral decomposition of EHM, in section \ref{sec:EHM}. In region ${\cal R}_1$, EHM has positive Lyapunov for all $\alpha$. Part (1) then follows from Theorem \ref{thm:dimspe=1} directly. Region ${\cal R_2}$ is actually where EHM has purely absolutely continuous measure for all $\alpha$ and a.e. $\theta$, see \cite{AJM} and Theorem \ref{thm:EHMspe} in section \ref{sec:EHM}. In view of Definition  (\ref{def:dimspe}), it is well known that if a measure is absolutely continuous w.r.t. Lebesgue measure, then it has full spectral dimension. Part (2) is then a direct consequence of a.c. spectrum and this fact.  We list part (2) here for completeness only. ${\cal R}_3$ is the region where EHM has zero Lyapunov exponent  and purely singular continuous spectrum for almost all $(\theta,\alpha)$,  part (3) follows from Theorem \ref{thm:speconti} and some technical improvements of the $(\Lambda,\beta)$ bound for analytic sampling functions. We will discuss in details about these three parts in section \ref{sec:EHM}. \\

Full spectral dimensionality is defined through the boundary behavior of the Borel
transform of the spectral measure. It implies a range of properties,
in particular, maximal packing dimension and quasiballistic quantum
dynamics. Thus our criterion links  way a purely analytic
property of the spectral measure to arithmetic property of the
frequency in a sharp.  In particular, consider $H_{\lambda,\alpha,\theta}$ , the extended Harper's model (EHM) given in (\ref{def:EHM}). In this part, we will focus on EHM and discuss the consequences of the full spectral dimensionality in terms of these explicit parameters. 

Recall that the Hausdorff/packing dimension of a (Borel) measure $\mu$, namely, ${\rm dim}_{\rm H}(\mu)/{\rm dim}_{\rm P}(\mu)$ is defined through the  $\limsup/\liminf$ ($\mu$ almost everywhere) of its $\gamma$-derivative $$\lim_{\varepsilon\downarrow 0}\frac{\ln \mu(E-\varepsilon,E+\varepsilon)}{\ln \varepsilon}.$$ 

If the $\liminf$ is replaced by $\limsup$ in the definition (\ref{def:dimspe}), we can define correspondingly the lower spectral
dimension ${\rm \underline {dim}}_{\rm spe}(\mu)$. It is well known (see e.g. \cite{DGY,Fal,JZ}) the relation between these fractal dimensions is 
${\rm dim}_{\rm H}(\mu)={\rm \underline {dim}}_{\rm spe}(\mu)\le {\rm {dim}}_{\rm spe}(\mu)\le {\rm dim}_{\rm P}(\mu)$.\footnote{In contrast to the Hausdorff
	dimension, the relation for the packing dimension only goes in one
	direction.} 
Therefore, lower bounds on spectral dimension lead to lower bounds on
packing dimension, thus also for the packing/upper box counting dimensions of the spectrum
as a set. We obtain corresponding non-trivial results for all the above quantities. The lower bounds also provide explicit examples where the spectral measure has different Hausdorff and packing dimension. 

Lower bounds on  spectral dimension also have immediate applications to
the lower bounds on quantum dynamics. Let $\delta_j\in\ell^2(\Z)$ be the delta vector in the usual sense. For $p > 0$, define
\begin{equation}\label{def:mom}
\langle|X|_{\delta_0}^p\rangle(T)=\frac{2}{T}\int_0^\infty e^{-2t/T}\sum_n|n|^p|\langle e^{- itH}\delta_0,\delta_n\rangle|^2.
\end{equation}
The power law of $\langle|X|_{\delta_0}^p\rangle(T)$ characterizes the propagation rate of $e^{- itH}\delta_0$. Define the upper/lower transport exponents to be
\begin{equation}\label{def:dyn-exp}
\beta_{\delta_0}^+(p)=\limsup_{T\to\infty}\frac{\ln\langle|X|_{\delta_0}^p\rangle(T)}{p\ln T}, \ \
\beta_{\delta_0}^-(p)=\liminf_{T\to\infty}\frac{\ln\langle|X|_{\delta_0}^p\rangle(T)}{p\ln T}.
\end{equation}
$\beta_{\delta_0}^-(p)=1$  for all $p>0$ corresponds to ballistic motion, $\beta_{\delta_0}^+(p)=1$  for all $p>0$  corresponds to quasiballistic motion.  $\beta_{\delta_0}^-(p)=0$ somtimes is called quasilocalized motion.  It was proved in \cite{GS99} that  $\beta_{\delta_0}^+(p)\ge {\rm dim}_{\rm P}(\mu), \ \forall p>0$. In view of  Corollary \ref{cor:EHM},  we have:
\begin{cor}\label{cor:EHM-2}
	Let $\mu_{\lambda,\alpha,\theta}$ be the spectral measure of EHM: $H_{\lambda,\alpha,\theta}$ defined in (\ref{def:EHM}). For any $\alpha\in[0,1]$, there is a full measure set $\Theta=\Theta(\alpha)\subset\T$ such that for any $\theta\in\Theta$,  $H_{\lambda,\alpha,\theta}$ has full packing dimension of $\mu_{\lambda,\alpha,\theta}$ and quasiballistic motion if  
	\begin{enumerate}
	\item  $\lambda\in {\cal R}_1$ and $\beta(\alpha)=\infty$.
	
	\item  $\lambda\in {\cal R}_2$ and for all $\alpha\in[0,1]$.
	
	\item  $\lambda\in {\cal R}_3$ and $\beta(\alpha)>0$.
	
\end{enumerate}
\end{cor} 

Hausdorff dimension of the spectral measure is always equal to
zero for a.e. phase for any ergodic operator \cite{S07} in
the regime of positive Lyapunov exponents. Combining the Lyapunov exponent formula of $H_{\lambda,\alpha,\theta}$ (see (\ref{def:EHMLyp})) with the result of Simon in \cite{S07},  we have ${\rm dim}_{\rm H}(\mu_{\lambda,\alpha,\theta})=0$ for $\lambda>1$,  a.e. $\theta$ and any $\alpha$.  In view of part (2) of Corollary \ref{cor:EHM-2}, for $\lambda>1$, a.e. $\theta$ and $\beta(\alpha)=\infty$, we have $0={\rm dim}_{\rm H}(\mu_{\lambda,\alpha,\theta})<{\rm dim}_{\rm P}(\mu_{\lambda,\alpha,\theta})=1.$ 


We are also interested in the fractal dimensional properties of the density states measure and the dimension of the spectrum as a set. Let ${\rm d}N_{\lambda,\alpha}$ be the density states measure and $\Sigma_{\lambda,\alpha}$ be the spectrum of $H_{\lambda,\alpha,\theta}$. For irrational $\alpha$, they are both $\theta$ independent. It is well known that 
\begin{align}
{\rm d}N_{\lambda,\alpha}=\E_\theta(\mu_{\lambda,\alpha,\theta})
\end{align}
and $\Sigma_{\lambda,\alpha}= {\rm supp_{top}}({\rm d}N_{\lambda,\alpha}).$ By these relations and the general  properties of the packing dimension of a measure and its topological support (see e.g. \cite{Fal}), Corollary \ref{cor:EHM-2} implies that 
\begin{align}
{\rm dim}_{\rm P}({\rm d}N_{\lambda,\alpha})={\rm dim}_{\rm P}(\Sigma_{\lambda,\alpha})=1
\end{align}
 in the corresponding parameter regions where EHM has full packing dimension.
 
For the dynamical transport part, Last in \cite{L1} proved that almost Mathieu operator with an appropriate Liouville frequency has quasiballistic motion for the first time. In general, quasiballistic property is a $G_\delta$ in any regular space, see e.g. \cite{S95,GKT}, thus this was known for (unspecified) topologically generic frequencies. In \cite{JZ}, the authors gave a precise
arithmetic condition on $\alpha$ for the quasiballistic motion depending on whether or not Lyapunov exponent vanishes
 in the quasiperiodic Schr\"odinger setting. Here, we prodive the parametric conditions for the EHM. The conclusions can also be extended directly to more general singular Jacobi operators with analytic quasiperiodic potentials. 

\section{Preliminaries}\label{sec:prelim}
We recall some commonly used notations for reader's convenience. We denote  ${ L}^{\infty}(\T,\R)$ and ${L}^{\infty}(\T,\C)$ to be the space of all 1-periodic bounded functions, taking values in $\R$ and $\C$ respectively. Denote the usual ${L}^{\infty}$ norm in both spaces by $	\|f\|_{\infty}:=\sup_{x\in\T}|f(x)|.$ Note we only require the diagonal potential function $v$ to be real valued functions, all the other sampling function are allowed to take value in $\C$. We do not emphasize the real/complex value anymore unless necessary. Denote ${L}^{1}(\T,\C)$ to be the usual Lebesgue space with the 1-norm $\|f\|_1:=\int_\T |f(\theta)|{\rm d}\theta$. Denote  $C^{\omega}(\T,\C)$ to be the space of all 1-periodic analytic functions and denote $C^{k}(\T,\C)$ to be the space of all  functions with continuous $k$-th order derivatives for all $k=0,1,\cdots,\infty$. 
 We denote ${\rm Lip}(\T,\C)$ to be the space of all 1-periodic Lipschitz continuous functions, induced with the Lipschitz norm given by:
\begin{align}\label{def:Lipnorm}
	\|f\|_{\rm Lip}:=\|f\|_\infty+\sup_{x,y\in\T}\frac{|f(x)-f(y)|}{|x-y|}.
\end{align}
We identify the sequence $u=\{u_n\}_{n\in\Z}$ with $u_n$ whenever it is clear that $n$ is the index. Denote the $\ell^{\infty}$ norm of $u\in \ell^{\infty}(\Z,\C)$ by $\|u\|_{\infty}:=\sup_{n\in\Z}|u_n|$. We will denote the distance on $\T^1$ by $\|\theta\|_{\T}:=\inf_{n\in\Z}|\theta-n|$ and may drop the subindex $\|\cdot\|=\|\cdot\|_{\T}$ whenever it is clear. 
\subsection{Transfer matrices and Lyapunov exponents}\label{sec:Transfer}

Let $H$ be given as in (\ref{def:Jacobi}):
\begin{equation}
(Hu)_n=w_n u_{n+1} +\overline{w}_{n-1} u_{n-1}+v_n u_n,\ n\in\Z.
\end{equation}
The eigenvalue equation $Hu=Eu$ can be rewritten via the following skew product:
 \begin{align} \label{eq:eigeneq}
\left (\begin{matrix}  u_{n+1} \\ u_{n} \end{matrix} \right )= A_n(E)\left (\begin{matrix}  u_{n} \\  u_{n-1}  \end{matrix} \right ),
 \end{align}
where
\begin{align} \label{def:AnDn}
A_n(E)=\frac{1}{w_n}D_n(E),\ D_n(E)=\left (\begin{matrix}  E-v_n   &  -\bar{w}_{n-1} \\ w_n  &   0 \end{matrix} \right ).
\end{align}
For $n\in\N^+$ and $m\in\Z$, define the n-step transfer matrix at  position $m$ to be
\begin{equation}\label{def:Anm}
A(n,m;E)=\prod_{j=m}^{n+m-1} A_j(E),
\end{equation}
\begin{equation}\label{def:Dnm}
D(n,m;E)=\prod_{j=m}^{n+m-1} D_j(E).
\end{equation}
We denote the scalar product of $w_n$ by the similar notation:
\begin{equation}\label{def:wnm}
w(n,m)=\prod_{j=m}^{n+m-1} w_j,\ \ m\in\Z,\ n\in\N^+
\end{equation}
and denote
\begin{align}
&A(n;E)=A(n,1;E), n>0;\  A(0;E)=Id;\ A(n;E)=A^{-1}(-n,n+1;E), n<0, \label{def:An}\\
&D(n;E)=D(n,1;E), n>0;\  D(0;E)=Id;\ D(n;E)=D^{-1}(-n,n+1;E), n<0,\label{def:Dn}
\end{align}
\begin{equation}\label{def:cn}
c(n)=c(n,1), n>0,
\end{equation}
for simplicity.\\

The (upper) Lyapunov exponent characterizes the grow(decay) rate of the norm of the transfer matrix $\|A(n,m)\|$, it will be convenient to introduce the Lyapunov exponent by
using the dynamical notations. We refer readers to \cite{JZ,AJ11} and references therein for the general definition of the Lyapunov exponent of linear skew product. In this part, we will restrict ourselves to the quasiperiodic cocycles. Let $\alpha\in\R\backslash\Q$ and $A:\T\mapsto {\rm GL}(2,\mathbb{C})$. We call $(\alpha,A)$ a (complex) cocycle. 
 In view of (\ref{def:Anm}) and (\ref{def:An}), denote the transfer matrix in the quasiperiodic cocycle case by
 \begin{equation}\label{def:Anmtheta}
 A(n;\theta,\alpha)=\prod_{j=1}^{n} A\big(\theta+(j-1)\alpha\big),\ \theta\in\T, n\in\N^+.
 \end{equation}
 
The
Lyapunov exponent is given by the formula:
\begin{equation}\label{def:Lyp}
L(A,\alpha)=\lim_{n\to +\infty}\frac{1}{n}\int_\T\ln\| A(n;\theta,\alpha)\|{\rm d}\theta=
\inf_{n>0}\frac{1}{n}\int_\T\ln\| A(n;\theta,\alpha)\|{\rm d}\theta.
\end{equation}
For irrational $\alpha$, the point-wise limit $L(A,\alpha)=\lim_{n\to +\infty}\frac{1}{n}\ln\| A(n;\theta,\alpha)\|$ also hold true for a.e. $\theta\in\T$ by subadditive ergodic theory. 

By uniquely ergodicity of the irrational rotations we have the following uniform upper bound (in $\theta$) for both matrix and scalar cases:
\begin{lem}[e.g. \cite{Fur,JMv2}]\label{lem:Furman}
If $A\in C^0(\T,{\rm GL}(2,\mathbb{C}))$, then 
\begin{align}\label{eq:limsup-A}
\limsup_{n\to +\infty}\frac{1}{n}\ln\| A(n;\theta,\alpha)\|\le L(A,\alpha)
\end{align}
uniformly in $\theta\in\T$. 

If $a\in C^0(\T,\mathbb{C})$ and $\ln|a(\theta)|\in L^1(\T)$, then 
\begin{align}\label{eq:limsup-a}
\limsup_{n\to +\infty}\frac{1}{n}\ln|\prod_{j=1}^{n} a\big(\theta+(j-1)\alpha\big)|\le \int_\T\ln|a(\theta)|{\rm d}\theta\end{align}
uniformly in $\theta\in\T$. 
\end{lem}
\begin{rem}
	If $a\in C^0(\T,\mathbb{C})$ has no zeros, then $\frac{1}{a(\theta)}$ is also continuous. By (\ref{eq:limsup-a}), we have 
	\begin{align}\label{eq:limsup-a-1}
	&\frac{1}{n}\ln\Big|\prod_{j=1}^{n} \frac{1}{a\big(\theta+(j-1)\alpha\big)}\Big|\le \int_\T\ln\Big|\frac{1}{a(\theta)}\Big|{\rm d}\theta+\epsilon\\
	\Longleftrightarrow &
	\prod_{j=1}^{n} |a\big(\theta+(j-1)\alpha\big)|\ge e^{n\big(\int_\T\ln |a(\theta)|{\rm d}\theta\,-\epsilon\big)}\end{align}
	for $n>n_0(\epsilon)$ (uniform in $\theta$).  This immediately gives the desired lower bound in (\ref{eq:Lambda-bound}) in a uniform way. If  $a(\theta)$ has zeros, there is no such uniform lower bound for the scalar product anymore. One technical achievement in the paper is, with some mild assumptions on the non-degeneracy of  the zeros,  we are able to get a weakened version of (\ref{eq:limsup-a-1}) (see Lemma \ref{lem:LipBound}), which will be sufficient for the spectral continuity.
		
\end{rem}

\subsection{The Weyl-Titchmarsh $m$-function and subordinacy theory}\label{sec:mMmu}
The boundary behavior of the $m$ function is linked to the power law of the half line solution and the growth of the transfer matrix norm $A(n,m;E)$ via the well known Gilbert-Pearson subordinacy theory \cite{G89,GP87}. 
We give a brief review on $m$-function and the subordinacy theory. More details can be found, e.g., in \cite{CL90}.

Let $H$ be as in (\ref{def:Jacobi}) and
$z=E+ i\varepsilon\in\mathbb{C}.$ Consider equation
\begin{equation}\label{eq:Hu=zu}
    Hu=zu.
\end{equation}
with the family of normalized phase boundary conditions:
\begin{equation}\label{eq:u-half-sln}
    u^{\varphi}_0\cos\varphi+u^{\varphi}_1\sin\varphi=0,\ -\pi/2<\varphi<\pi/2,\  |u^{\varphi}_0|^2+|u^{\varphi}_1|^2=1.
\end{equation}
Let $\Z^+=\{1,2,3\cdots\}$ and $Z^-=\{\cdots,-2,-1,0\}$. Denote by $u^{\varphi}=\{u^{\varphi}_j\}_{j\ge0}$ the right half line solution on $\Z^+$ of (\ref{eq:Hu=zu}) with  boundary condition (\ref{eq:u-half-sln}) and by $u^{\varphi,-}=\{u^{\varphi,-}_j\}_{j\le0}$ the left half line solution on $\Z^-$ of the same equation. Also denote by $v^{\varphi}$ and $v^{\varphi,-}$  the right and left half line solutions of (\ref{eq:Hu=zu}) with the orthogonal boundary
conditions to $u^{\varphi}$ and $u^{\varphi,-}$, i.e., $v^{\varphi}=u^{\varphi+\pi/2}$,$v^{\varphi,-}=u^{\varphi+\pi/2,-}$. For any function $u:\Z^+\to\mathbb{C}$ we denote by $\|u\|_{\ell}$ the norm of $u$ over a lattice
interval of length $\ell$; that is
\begin{equation}\label{eq:ul}
\|u\|_{\ell}=\Big[\sum_{n=1}^{[\ell]}|u_n|^2+(\ell-[\ell])|u_{[\ell]+1}|^2\Big]^{1/2}.
\end{equation}
Similarly, for $u:\Z^-\to\mathbb{C}$, we define
\begin{equation}\label{eq:ul-}
\|u\|_{\ell}=\Big[\sum_{n=1}^{[\ell]-1}|u_{-n}|^2+(\ell-[\ell])|u_{-[\ell]}|^2\Big]^{1/2}.
\end{equation}
For any $\varepsilon>0$, let $\ell=\ell(\varphi,\varepsilon,E)$ be
\begin{equation}\label{eq:subor-l}
    \|u^{\varphi}\|_{\ell(\varphi,\varepsilon)}\|v^{\varphi}\|_{\ell(\varphi,\varepsilon)}=\frac{1}{2\varepsilon}.
\end{equation}
$\ell^-(\varphi)$ is defined through the same equation by $u^{\varphi,-},v^{\varphi,-}$. It is easy to check
\begin{equation}\label{eq:uv0}
   \|u^{\varphi}\|_{\ell}\cdot \|v^{\varphi}\|_{\ell}\ge \frac{1}{2}([\ell]-1).
\end{equation}

Let $m_{\varphi}(z):\mathbb{C}^+\mapsto\mathbb{C}^+$ and
$m_{\varphi}^-(z):\mathbb{C}^+\mapsto\mathbb{C}^+$ the right and left Weyl-Titchmarsh
m-functions (half line) associated with the boundary condition (\ref{eq:u-half-sln}). Let
$m=m_0$ and $m^-=m^-_0$ be  the half line m-functions corresponding to the Dirichlet boundary conditions.
The following quantitative subordinate theory was proved in \cite{JL1}, well known as Jitomirskaya-Last inequality.  
\begin{lemma}[Jitomirskaya-Last inequality, Theorem 1.1 in \cite{JL1}] \label{lem:JLineq}
For $E\in\R$ and $\varepsilon>0$, the following inequality holds for any $\varphi\in(-\frac{\pi}{2},\frac{\pi}{2}]$:
\begin{equation}\label{eq:jl}
    \frac{5-\sqrt{24}}{|m_{\varphi}(E+i\varepsilon)|}
<\frac{\|u^{\varphi}\|_{\ell(\varphi,\varepsilon)}}{\|v^{\varphi}\|_{\ell(\varphi,\varepsilon)}}
<\frac{5+\sqrt{24}}{|m_{\varphi}(E+i\varepsilon)|}.
\end{equation}

\end{lemma}

There is also one general statement about the existence of generalized eigenfunctions with sub-linear growth in its $\ell$-norm:
\begin{lemma}[\cite{LS99}]\label{lem:ls99}
	For $\mu_\theta$-a.e. $E$, there exists $\varphi\in (-\pi/2,\pi/2]$ such that $u^{\varphi}$ and $u^{\varphi,-}$ both obey
	\begin{equation}\label{eq:ls99}
	\limsup_{\ell\to\infty}\frac{\|u\|_{\ell}}{\ell^{1/2}\ln \ell}<\infty.
	\end{equation}
\end{lemma}
This inequality provides us an upper bound for the $\ell$-norm of the solution, which is crucial in the proof of the spectral singularity. \\

The next proposition relates the whole line m-function $M$ and half line m-function $m_\varphi$, which can be found 
in \cite{DKL00}.

\begin{prop}[Corollary 21 in \cite{DKL00}] \label{prop:DKL}
Fix $E\in\R$ and $\varepsilon>0$,
\begin{equation}\label{eq:DKL}
    |M(E+ i\varepsilon)|\le\sup_{\varphi}|m_{\varphi}(E+ i\varepsilon)|.
\end{equation}
\end{prop}
By this proposition, to bound $M$ from above and get spectral continuity as in (\ref{def:speconti}), it is enough to obtain uniform upper bounds of $m_{\varphi}$ in boundary condition $\varphi$ for the right half line problem.\\

For spectral singularity, we need to consider both $m_{\varphi}(z)$ and $m_{\varphi}^-(z)$. Let $(U\psi)_n=\psi_{-n+1},\ n\in\Z$ be a a unitary operator on ${\ell}^2(\Z)$. Let $\widetilde{H}=UHU^{-1}$. Denote by
$\widetilde{m},\widetilde{m}_\varphi,\widetilde{u}^{\varphi}$ and
$\widetilde{\ell}(\varphi)$, correspondingly,
${m},{m}_\varphi,{u}^{\varphi}$ and ${\ell}(\varphi)$ of the operator
$\widetilde{H}$. The following facts are well known in the past literatures(see e.g. section 3, \cite{JL2}). For any $\varphi\in(-\pi/2,\pi/2]$,
\begin{equation}\label{eq:Mm+m-}
    M(z)=\frac{m_\varphi(z)\widetilde{m}_{\pi/2-\varphi}-1}{m_\varphi(z)+\widetilde{m}_{\pi/2-\varphi}}
\end{equation}
and
\begin{equation}\label{eq:L+-}
    \widetilde{\ell}(\pi/2-\varphi)={\ell}^-(\varphi),\  \ \|u\|_{\ell}=\|Uu\|_{\ell}.
\end{equation}

In view of \eqref{def:speSingular}, a direct consequence of (\ref{eq:Mm+m-}) is (e.g. Lemma 5 in \cite{JL2}):
\begin{lemma}\label{lem:Mmm}
For any $0<\gamma<1$, suppose that there exists a $\varphi\in(-\pi/2,\pi/2]$ such that for $\mu$-a.e. $E$ in some Borel set $S$, we have that $\liminf_{\varepsilon\to0} \varepsilon^{1-\gamma}|m_\varphi(E+ i\varepsilon)|=\infty$ and $\liminf_{\varepsilon\to0} \varepsilon^{1-\gamma}|\widetilde{m}_{\pi/2-\varphi}(E+ i\varepsilon)|=\infty$.
Then for $\mu$-a.e. $E$ in $S$, $\liminf_{\varepsilon\to0} \varepsilon^{1-\gamma}|M(E+ i\varepsilon)|=\infty$, namely,  the restriction $\mu(S\cap\cdot)$ is $\gamma$-spectral singular.
\end{lemma}

\subsection{Continued fraction}\label{sec:beta}
An important tool in the study of quasiperiodic sequence is the continued fraction expansion of irrational numbers. Let $\alpha\in \T\setminus \Q$, $\alpha$ has the following unique expression with $a_n\in \N$:
\begin{align}\label{def:continuefractionan}
\alpha=\frac{1}{a_1+\frac{1}{a_2+\frac{1}{a_3+\cdots}}}.
\end{align}
Let
\begin{align}\label{def:pnqn}
\frac{p_n}{q_n}=\frac{1}{a_1+\frac{1}{a_2+\frac{1}{\cdots+\frac{1}{a_n}}}}
\end{align}
be the continued fraction approximants of $\alpha$.
Let
\begin{align*}
\beta(\alpha)=\limsup_{n\rightarrow\infty}\frac{\ln{q_{n+1}}}{q_n}.
\end{align*}
$\beta(\alpha)$ being large means $\alpha$ can be approximated very well by a sequence of rational numbers.
Let us mention that $\{\alpha: \beta(\alpha)=0\}$ is a full measure set.

The following properties about continued fraction expansion are well-known:
\begin{align}\label{eq:qnqn+1}
\frac{1}{2q_{n+1}}\leq \|q_n\alpha\|_{\T}\leq \frac{1}{q_{n+1}}.
\end{align}
For any $q_n\leq |k|<q_{n+1}$,
\begin{align}\label{eq:qnkqn+1}
\|q_n\alpha\|_{\T}\leq \|k\alpha\|_{\T}.
\end{align}
Combining definition of $\beta(\alpha)$ (\ref{def:beta}) with (\ref{eq:qnkqn+1}), we have:
If $\beta(\alpha)=0$, then for any $\delta>0$, for $|k|$ large, the following inequality holds:
\begin{align}\label{eq:beta=0epsilon}
\|k\alpha\|_{\T}>e^{-\delta |k|}.
\end{align}

\subsection{More about $\beta$-almost periodicity and the $(\Lambda,\beta)$ bound.}\label{sec:MoreLambdaBound}
In this paper, we consider bounded  sequences $v_n,w_n$ (for example, $v_n,w_n$ are both generated by some smooth sampling functions). Let $D(n,m)$ and $w(n,m)$ be defined as in (\ref{def:Dnm}) and (\ref{def:wnm}). The mild assumption on $v_n,w_n$ yields the following trivial upper bound for $D(n,m)$ and $w(n,m)$: there is $\Lambda_0=\Lambda_0(\|v\|_\infty,\|w\|_\infty)>0$ such that
  for any $n\in\N$, and any $E\in{\cal N}:={\cal N}(H)$,
\begin{align}
\sup_{m\in\Z} \|D(n,m;E)\| \leq e^{\Lambda_0 n}, \label{def:D-upper}\\
\sup_{m\in\Z} |w(n,m)| \leq e^{\Lambda_0 n}. \label{def:wnm-upper}
\end{align}
Suppose $w_n$ has $(\Lambda,\beta)$-$q$ bound as in (\ref{eq:Lambda-bound}). Without of generality, we assume $\Lambda_0=\Lambda$ for simplicity. 
For $1\le r\le q$ and $m\in\Z$, write $w(q,m)={w(q-r,m+r)}w(r,m)$. Combine (\ref{eq:Lambda-bound}) with the upper bound (\ref{def:wnm-upper}), we have
 \begin{align}\label{eq:wrm-lower}
\min_{|m|\le e^{\delta\beta\,q}} |w(r,m)| \ge e^{-2\Lambda\,q}, \ \ 1\le r< q.
\end{align}
In particular, $r=1$ gives
 \begin{align}\label{eq:wm-lower}
\min_{|m|\le e^{\delta\beta\,q}} |w_m| \ge e^{-2\Lambda\,q}.
\end{align}
Assume further $w_n$ has $\beta$-$q$ almost periodicity as in (\ref{eq:a-beta}), by (\ref{eq:wm-lower}),  $\beta$-$q$ periodicity can be strengthened as, 
 \begin{align}\label{eq:wq/wq}
\max_{|m|\le e^{\delta\beta\,q}}\Big|\frac{w_{m\pm q}}{w_{m}}-1\Big| < e^{-(\beta-2\Lambda) q}
\end{align}

We also abuse the notation frequently by saying
the operator $H$ or the transfer matrix $A(n,E)$ has  $\beta$-almost periodicity and $(\Lambda,\beta)$ boundedness
if the corresponding $v_n,w_n$ has $\beta$  almost periodicity and $(\Lambda,\beta)$ boundedness .

The lower bound on $w(n,m)$ and upper bound on $D(n,m)$ also imply that for any $E\in{\cal N}$, and $|m|\le e^{\delta\beta\,q}$
 \begin{align}\label{eq:Arm-upper}
\|A(q,m)\|< e^{2\Lambda q},\ \ \max_{0\le r< q}\|A(r,m)\|< e^{3\Lambda q}. 
\end{align}

Assume now $\beta$-almost periodicity and $(\Lambda,\beta)$ bound hold true for the sequence $q_n\to\infty$, we will use these induced bounds (\ref{eq:wrm-lower})-(\ref{eq:Arm-upper}) for the sequence $q_n$ frequently. 


\section{$(\Lambda,\beta)$ bound for quasiperiodic smooth sequence and the proof of Corollary \ref{cor:LipJacobi}}\label{sec:LambdaBound}
Assume we have a the quasiperiodic sequence $v(\theta+n\alpha)$ generated by a Lipschitz sampling function $v$. Let $q_n$ be given as in (\ref{def:beta}). By (\ref{def:beta}), for any $0<\beta<\beta(\alpha)/2$, there is a subsequence $q_{n_k}$ such that $\ln q_{n_k+1}>2\beta q_{n_k}$. 
Then for any $\theta,j$ and $1\le n\le q_{n_k}$,
 \begin{align}\label{eq:Lip-beta}
|v\big(\theta+m\alpha\big)-v\big(\theta+(m\pm q_{n_k})\alpha\big)|\le \|v\|_{\rm Lip}\cdot\|q_{n_k}\alpha\|\le \|v\|_{\rm Lip}\cdot\frac{1}{q_{n_k+1}}\le \|v\|_{\rm Lip}\cdot e^{-2\beta
	q_{n_k}}\le e^{-\beta
	q_{n_k}},
\end{align}
provided $q_{n_k}$ large. Same computation works for $c$. Therefore, $v(\theta+n\alpha)$ and $c(\theta+n\alpha)$ are $\beta$-almost periodic for Lipschitz continuous $v,c$.

The more challenging part is the $(\Lambda,\beta)$ bound on $c(\theta+n\alpha)$, where we need some further assumption on $c$. We will focus on this throughout the rest of this section.

The key ingredient for the proof of the $(\Lambda,\beta)$ bound is the following lemma in \cite{AJ09}:
  \begin{lemma} \label{lem:AJ09}
 Let $\alpha\in \R\backslash\Q $,\ $\theta\in\R$ and $0\leq j_0 \leq q_{n}-1$ be such that 
 $$\mid \sin \pi(\theta+j_{0}\alpha)\mid = \inf_{0\leq j \leq q_{n}-1} \mid \sin \pi(\theta+j\alpha)\mid ,$$
 then for some absolute constant $C>0$,
 $$-C\ln q_{n} \leq \sum_{j=0,j\neq j_0}^{q_{n}-1} \ln \mid \sin \pi (\theta+j\alpha) \mid+(q_{n}-1)\ln2 \leq C\ln q_n. $$
 \end{lemma}
This lemma was used  in \cite{JY} to  prove some optimal singular continuous spectrum results. By extending the argument in \cite{JY} to exponentially many periods, we are able to prove  $(\Lambda,\beta)$ bound for any analytic sampling function. Actually, we can deal with more general sampling functions with much weaker regularities. Define
\begin{align}\label{def:PL}
{\cal {F}}(\T,\mathbb{C}):=\Big\{&c\in {L}^{\infty}({\T,\mathbb{C}}):\ \exists m \in\N^+, \theta_\ell\in\T, \tau_\ell\in(0,1],\ell=1,\cdots,m \nonumber
\\
&\qquad \textrm{ such that}\ \ 
g(\theta):=\frac{c(\theta)}{\prod_{\ell=1}^m|\sin \pi(\theta-\theta_{\ell})|^{\tau_{\ell}}}
\in {L}^{\infty}(\T,\mathbb{C})
\ {\rm and}\ \inf_\T |g(\theta)|>0. 
\Big\} 
\end{align}

Suppose $c(\theta)\in {\cal {F}}(\T,\mathbb{C})$ with $\theta_\ell$ and $g(\theta)$ given as in (\ref{def:PL}) such that
\begin{align}
c(\theta)=g(\theta)\,\prod_{\ell=1}^m|\sin\pi(\theta-\theta_{\ell})|^{\tau_{\ell}}.
\end{align}
Clearly, $\ln |g(\theta)|\in L^1(\T)$. By the well known integral $\int_\T\ln |\sin\pi\theta\,|{\rm d}\theta=-\ln2$, it is easy to check that  $\ln |c(\theta)|\in L^1(\T)$ and is linked to $\ln|g|$ by:
\begin{align}
\label{eq:int-cg}\int_\T\ln |c(\theta)|{\rm d}\theta=\int_\T\ln |g(\theta)|{\rm d}\theta-\ln2 \sum_{\ell=1}^m\tau_\ell.
\end{align}

The following technique lemma shows that any sampling function in  ${\cal {F}}(\T,\mathbb{C})$ with an irrational force can generate a $(\Lambda,\beta)$ bounded sequence.
\begin{lemma}\label{lem:LipBound}
Assume that there exists $m \in\N^+, \theta_{\ell}\in\T, \tau_{\ell}\in(0,1],\ell=1,\cdots,m,$  $g(\theta)\in {L}^{\infty}(\T,\mathbb{C})$  such that $\inf_\T |g(\theta)|>0$ and 
\begin{align}\label{eq:cinPL}
c(\theta)=g(\theta)\prod_{{\ell}=1}^m|\sin \pi(\theta-\theta_{\ell})|^{\tau_{\ell}}.
\end{align}
Then for any $\alpha$ with  $0<2\beta<\beta(\alpha)$ and $0<\delta<\frac{2\sum_{\ell=1}^m\tau_{\ell}}{1+\sum_{\ell=1}^m\tau_{\ell}}$, there is a sequence $q_n\to \infty$  and a full Lebesgue measure set $\Theta=\Theta(\alpha,\theta_1,\cdots,\theta_m)$ such that 
for any $\theta\in\Theta$ and $q_n$ large enough\footnote{The sequence itself only depends $\beta(\alpha)$, while the largeness depends on $\theta,\alpha,\beta,\delta,\tau$.},  $c(\theta+n\alpha)$ satisfies:
\begin{align}\label{eq:Lambda-bound-Lip}
\min_{|k|\le q_{n}^{-1}e^{\delta\beta\,q_n}}\prod_{j=kq_n}^{(k+1)q_n-1}  |c(\theta+j\alpha)|>e^{-\Lambda_1\,q_n},
\end{align}
where 
\begin{align}\label{eq:tempLambda1}
\Lambda_1:=\Lambda_1(\tau,g,\delta\beta)=
\ln2\sum_{\ell=1}^m\tau_{\ell}-\ln\big(\inf_\T|g(\theta)|\big)+\delta^2 \min\{\beta,1\}.
\end{align}

Assume further $g(\theta)\in C^0({\T,\mathbb{C}})$ and $\ln |g(\theta)|\in L^1(\T)$,  $\Lambda_1$ 
in (\ref{eq:tempLambda1}) can be replaced by 
\begin{align}\label{eq:tempLambda1-2}
\Lambda_1=-\int_\T\ln |c(\theta)|{\rm d}\theta+2\delta^2 \min\{\beta,1\}.
\end{align}
Moreover, 
  $c(\theta+n\alpha)$ is $(\Lambda,\beta)$ bounded as defined in (\ref{eq:Lambda-bound}) such that 
  \begin{align}\label{eq:cmqn}
  \min_{|m|\le e^{\delta\beta\,q_n}}\prod_{j=m}^{m+q_n-1}   |c(\theta+j\alpha)|>e^{-\Lambda\,q_n},\ \ \theta\in\Theta
  \end{align}
  where 
\begin{align}
\Lambda=-\int_\T\ln |c(\theta)|^{}{\rm d}\theta+6\delta^2 \min\{\beta,1\} .\label{eq:tempLambda} 
\end{align}  
  
\end{lemma}
\begin{rem}
The above $\Lambda_1$ and $\Lambda$ can be negative in general, which makes (\ref{eq:Lambda-bound-Lip} and (\ref{eq:cmqn}) actually exponentially grow (instead of decay). This is natural since  there is actually a `large' scaling of size $\int \ln |g|\sim \int \ln |c|$ for the product in these cases. We are more interested in the case where $\int \ln |c|\le 0$ where $\Lambda_1$ and $\Lambda$ are indeed positive. In particular, it is always possible to re-scale $c(\theta)$ to make the logarithm average zero. This will lead to an arbitrarily small $\Lambda$ (positive) in (\ref{eq:tempLambda}). Combine this with the uniform Lyapunov upper bound for (\ref{def:D-upper}) (\cite{Fur}), we can get some refined results in the zero Lyapunov regime, e.g., the critical EHM model, about the spectral continuity and quasi-ballistic motion. See more discussion in the next Corollary \ref{cor:c-analytic-bound} and section \ref{sec:EHM} about Corollary \ref{cor:EHM}. 
\end{rem}

\begin{rem}
	It is an easy exercise that if $c(\theta)$ is $C^{k}$ continuous on $\T$ with finitely many zeros with non-degenerate $k$-th order derivatives, then $c(\theta)\in {\cal {F}}(\T,\mathbb{C})\bigcap {\rm Lip}(\T,\C)$ with a continuous $g(\theta)$ for all $k\ge 1$. By Lemma \ref{lem:LipBound}, there is $\Lambda=\Lambda(c)$ such that $c(\theta+n\alpha)$ is $(\Lambda,\beta)$ bounded for any $0<2\beta<\beta(\alpha)$ and a.e. $\theta\in\T$. From the proof of Lemma \ref{lem:LipBound}, see \eqref{eq:temp-qn},  the subsequence $q_{n_k}$ can be taken to be same as in the $\beta$-almost periodicity (\ref{eq:Lip-beta}). Therefore, Corollary \ref{cor:LipJacobi} follows from Theorem \ref{thm:speconti}. We omit the details here. An interesting question is whether the non-degenerate condition on $c$ can be weakened and what is the appropriate `non-degenerate' condition on any Lipschitz function such that (\ref{def:PL}) holds. 
\end{rem}
\proof 
Let $0<2\beta<\beta(\alpha)$ and $q_n$ be defined as in (\ref{def:beta}). For any $\delta>0$, let $q_{n_k}$ be the subsequence such that ${\ln q_{n_k+1}}\ge 2\beta\,{q_{n_k}}$. For simplicity, drop the subindex $n_k$ and denote the  subsequence still by $q_n$, i.e., 
\begin{align}\label{eq:temp-qn}
q_{n+1}>e^{2\beta\,q_n}
\end{align}
We also write $\tilde \beta=\min\{\beta,1\}$. It is obvious that for any $m\in\Z$ and $\theta\in\T$,
\begin{align}\label{eq:bound-g}
 \prod_{j=kq_n}^{(k+1)q_n-1}   |g(\theta+j\alpha)|>\big(\inf_\T|g(\theta)|\big)^{q_n}=e^{q_n\ln(\inf_\T|g(\theta)|)}.   
\end{align}
In view of (\ref{eq:cinPL}), it is enough to study the lower bound for each $\|\theta-\theta_{\ell}\|_\T$. For any $\alpha\in[0,1]\backslash\Q$ and any $\theta_\ell,\ell=1,\cdots,m$, let 
\begin{align}\label{eq:Theta_l}
\Theta_{\ell}:=\bigcup_{\gamma>0}\Big\{\theta\in\T:\, \|\theta-\theta_{\ell}+n\alpha\|_\T\ge \gamma |n|^{-2},\ \forall n\in\Z\backslash \{0\}\Big\}
\end{align}
It is well known that $\Theta_{\ell}$ is a full measure set. 
Let 
\begin{align}\label{eq:Theta}
\Theta:=\bigcap_{{\ell}=1}^m\Theta_{\ell}.
\end{align}
For any $\theta\in\Theta$ and  $1\le{\ell}\le m$, there is $\gamma_{\ell}=\gamma(\theta,\theta_{\ell},\alpha)>0$ such that 
\begin{align}
 \|\theta-\theta_{\ell}+n\alpha\|_\T\ge \frac{\gamma_{\ell}} {|n|^{2}},\ \forall n\in\Z\backslash \{0\}. \label{eq:thetaDC}
\end{align}

For all $1\le \ell  \le m$ and  $|k|<q_n^{-1}e^{\delta\beta q_n}$, 
let $j_{{\ell},k}\in [0,q_n)$ be such that the following holds:
$$|\sin\pi\big(\theta-\theta_{\ell}+kq_n\alpha+j_{\ell,m} \alpha \big)| 
=\inf_{0\leq j < q_n}|\sin\pi\big(\theta-\theta_l+kq_n\alpha+j\alpha\big)|.$$
By (\ref{eq:thetaDC}), for all $ j_{{\ell},k}$, 
\begin{align}
 \|\theta-\theta_{\ell}+j_{{\ell},k}\alpha\|_\T\ge \frac{\gamma_{\ell}} {|j_{{\ell},k}|^{2}}\ge \frac{\gamma_{\ell}} {q_n^{2}}.
\end{align}

Let $
\tau=\sum_{\ell=1}^m\tau_{\ell}
$. For all $|k|<q_n^{-1}e^{\delta\beta\,q_n}<e^{\delta\beta\,q_n}$, 
 we have that 
\begin{align}
\big|\sin\pi(\theta-\theta_{\ell}+kq_n\alpha+j_{{\ell},k}\alpha)\big|\ge
 \|\theta-\theta_{\ell}+kq_n\alpha+j_{{\ell},k}\alpha\|_\T &\ge  \|\theta-\theta_{\ell}+j_{{\ell},k}\alpha\|_\T-\|kq_n\alpha\|_\T  \nonumber\\
 &\ge \frac{\gamma_{\ell}} {q_n^{2}}-|k|\frac{1}{q_{n+1}}\nonumber \\
 &\ge \frac{\gamma_{\ell}} {q_n^{2}}-e^{\delta\beta q_n}{e^{-2\beta q_n}} \nonumber \\
 &\ge 2e^{-\tau^{-1}\delta^2\widetilde \beta q_n}-e^{-(2-\delta)\widetilde\beta q_n}\nonumber\\
 & \ge e^{-\tau^{-1}\delta^2\widetilde \beta q_n} \label{eq:tempTau}
\end{align}
provided $q_n^{-2}e^{\tau^{-1}\delta^2\tilde{\beta} q_n}\ge 2\gamma_{\ell}^{-1}$ and $2-\delta>\tau^{-1}\delta>\tau^{-1}\delta^2$. The latter gives the restriction on $\delta$ such that $\delta<\frac{2}{1+\tau^{-1}}$.

By Lemma \ref{lem:AJ09}, 
\begin{align}
 \prod_{j=0,j\neq j_{l,k}}^{q_n-1}|\sin\pi(\theta-\theta_{\ell}+kq_n\alpha+j\alpha)|
\ge e^{-(q_n-1)\ln2-C\ln q_n} \ge e^{-q_n\ln2-\tau^{-1}\delta^2\widetilde \beta q_n} \label{eq:tempTau2}
\end{align}
provided $C\ln q_n<\tau^{-1}\delta^2\widetilde \beta q_n$ where $C$ is the absolute constant in Lemma \ref{lem:AJ09} and $\tau$ is the same in (\ref{eq:tempTau}).

Now putting \eqref{eq:tempTau} and \eqref{eq:tempTau2} together, we have
\begin{align*}
& \prod_{j=kq_n}^{(k+1)q_n-1}\Big(\prod_{{\ell}=1}^m|\sin\pi(\theta+j\alpha-\theta_{\ell})|^{\tau_{\ell}}\Big)\\
 &= \prod_{{\ell}=1}^m\big(\prod_{j=0}^{q_n-1}|\sin\pi(\theta-\theta_{\ell}+kq_n\alpha+j\alpha)|\big)^{\tau_{\ell}}\\
 &=  \Big(\prod_{{\ell}=1}^m\big(\prod_{j=0,j\neq j_{l,k}}^{q_n-1}|\sin\pi(\theta-\theta_{\ell}+kq_n\alpha+j\alpha)|\big)^{\tau_{\ell}}\Big)
 \cdot \Big(\prod_{{\ell}=1}^m|\sin\pi(\theta-\theta_{\ell}+kq_n\alpha+j_{l,k}\alpha)|^{\tau_{\ell}}\Big)\\
 &\ge  \Big(\prod_{{\ell}=1}^m\big(e^{-q_n\ln2-\tau^{-1}\delta^2\widetilde \beta q_n}\big)^{\tau_{\ell}}\Big)
 \cdot \Big(\prod_{{\ell}=1}^m\big(e^{-\tau^{-1}\delta^2\widetilde \beta q_n}\big)^{\tau_{\ell}}\Big) \\
 &= e^{-q_n(\ln2\sum_{\ell=1}^m\tau_{\ell})-\delta^2\widetilde \beta q_n}.
\end{align*}

Combined with (\ref{eq:bound-g}), we have that for all ${|k|\le q_{n}^{-1}e^{\delta \beta\,q_n}}$,
\begin{align}\label{eq:temp-1}
\prod_{j=kq_n}^{(k+1)q_n-1}  |c(\theta+j\alpha)|
>e^{-\big(\ln2\sum_{\ell=1}^m\tau_{\ell}-\ln\inf_\T|g(\theta)|+ \delta^2\widetilde \beta\big)q_n}
\end{align}
provided $q_n>\tilde q=\tilde q\big(\max_{\ell}\gamma^{-1}_{\ell},\delta,\alpha,\sum_{\ell=1}^m\tau_{\ell}\big)$.\\

Assume further $g(\theta),c(\theta)\in C^0(\T,\mathbb{C})$. Since $\inf |g(\theta)|>0$, $\ln {|g(\theta)|^{-1}}$ is also continuous. By Lemma \ref{lem:Furman}, there is $n_0=n_0(\delta^2\widetilde \beta)$ such that the following upper bound holds uniform in $\theta\in\T$ for $n>n_0$:
\begin{align}
\frac{1}{n}\sum_{j=1}^{n}  \ln |g(\theta+j\alpha)|^{-1}\le \int_\T\ln |g(\theta)|^{-1}{\rm d}\theta+\delta^2\widetilde \beta.
\end{align}
In particular, for all $q_n\ge n_0$ and any $k\in\Z$ we have
\begin{align*}
\Big(\prod_{j=0}^{q_n-1}|g(\theta+kq_n\alpha+j\alpha)|^{-1}\Big)^{\frac{1}{q_n}}
\le e^{-\int_\T\ln |g(\theta)|{\rm d}\theta+\delta^2\widetilde \beta}
\Longrightarrow \prod_{j=kq_n}^{(k+1)q_n-1}|g(\theta+j\alpha)|\ge e^{q_n\big(\int_\T\ln |g(\theta)|{\rm d}\theta-\delta^2\widetilde \beta\big) }.
\end{align*}

Therefore, we can replace $\ln\big(\inf_\T|g(\theta)|\big)$ 
in (\ref{eq:temp-1}) by $\int_\T\ln |g(\theta)|{\rm d}\theta-\delta^2 \min\{\beta,1\}$. In view of (\ref{eq:int-cg}), we have
\begin{align}\label{eq:tempLambda1-3}
\Lambda_1=
\ln2\sum_{\ell=1}^m\tau_{\ell}-\int_\T\ln |g(\theta)|{\rm d}\theta+2\delta^2 \min\{\beta,1\}=-\int_\T\ln |c(\theta)|{\rm d}\theta+2\delta^2 \min\{\beta,1\},
\end{align}
which gives the desired expression of $\Lambda_1$ in (\ref{eq:tempLambda1-2}).

Let $\widetilde c(\theta)=c(\theta)\,e^{-\int_\T\ln|c(\theta)|{\rm d}\theta}$. It is easy to check that 
$\int_\T\ln|\widetilde c(\theta)|{\rm d}\theta=0$. By (\ref{eq:tempLambda1-3}), we have for all ${|k|\le q_{n}^{-1}e^{\delta \beta\,q_n}}$,
\begin{align}\label{eq:temp-3}
\prod_{j=kq_n}^{(k+1)q_n-1}  |\widetilde c(\theta+j\alpha)|
>e^{-2\delta^2\widetilde \beta\,q_n}.
\end{align}

By Lemma \ref{lem:Furman}, there is $r_0=r_0(\delta^2\widetilde \beta)\in\N$ such that for any $m\in\Z$ and $r\ge r_0$, 
\begin{align}
\prod_{j=m}^{m+r-1}|\widetilde c(\theta+j\alpha)|
\le e^{r\big(\int_\T\ln \widetilde c(\theta)|{\rm d}\theta+\delta^2\widetilde \beta\big)}
= e^{\delta^2\widetilde \beta\,r}.
\end{align}
For $0\le r<r_0$, we have the trivial upper bound $\prod_{j=m}^{m+r-1}  |\widetilde c(\theta+j\alpha)| \le e^{r\ln(\|\widetilde c\|_\infty)}\le e^{r_0\ln(\|\widetilde c\|_\infty+1)}$. Therefore, for any $m\in\Z$ and $1\le r \le q_n$,  
\begin{align}
\prod_{j=m}^{m+r-1}|\widetilde c(\theta+j\alpha)|
\le e^{\delta^2\widetilde \beta\,q_n}
\end{align}
provided $q_n\ge \delta^{-1}\widetilde\beta^{-1}r_0\ln(\|\widetilde c\|_\infty+1)$.

Then for any $|m|<e^{\delta \beta\,q_n}$, there is $k$ such that $kq_n\in (m,m+q_n]$, therefore,
\begin{align}
\prod_{j=m}^{m+q_n-1}  |\widetilde c(\theta+j\alpha)|&=\prod_{j=m}^{kq_n-1}  |\widetilde c(\theta+j\alpha)|\cdot \prod_{j=kq_n}^{m+q_n-1}  |\widetilde c(\theta+j\alpha)|\\
&=\frac{\prod_{j=(k-1)q_n}^{kq_n-1}  |\widetilde c(\theta+j\alpha)|} {\prod_{j=(k-1)q_n}^{m-1}  |\widetilde c(\theta+j\alpha)|}\cdot
\frac{\prod_{j=kq_n}^{(k+1)q_n-1}  |\widetilde c(\theta+j\alpha)|}{\prod_{j=m+q_n}^{(k+1)q_n-1}  |\widetilde c(\theta+j\alpha)|}\\
&>\frac {e^{-2\delta^2\widetilde \beta\,q_n}}{e^{\delta^2\widetilde \beta\,q_n}}\cdot \frac {e^{-2\delta^2\widetilde \beta\,q_n}}{e^{\delta^2\widetilde \beta\,q_n}} \\
&=e^{-6\delta^2\widetilde \beta\,q_n} .
\end{align}
Therefore, 
\begin{align}
\prod_{j=m}^{m+q_n-1}  |c(\theta+j\alpha)|&=e^{q_n\,\int_\T\ln|c(\theta)|{\rm d}\theta}\prod_{j=m}^{m+q_n-1}  |\widetilde c(\theta+j\alpha)|\nonumber \\
&\ge e^{-(-\int_\T\ln|c(\theta)|{\rm d}\theta+6\delta^2\widetilde \beta)\,q_n}=:
e^{-\Lambda\,q_n} \label{eq:tempLambda-comput},
\end{align}
as claimed.  
This completes the proof of Lemma \ref{lem:LipBound}.
\qed

As an explicit example, we have the following arbitrarily slow lower bound for the analytic case with zero $\ln$ mean. 
\begin{cor}\label{cor:c-analytic-bound}
Assume that $c(\theta)\in C^{\omega}(\T,\C)$ and $\int_\T\ln|c(\theta)|{\rm d}\theta=0$. Denote the all zeros
\footnote{Clearly, analytic function $c(\theta)$ only has finitely many zeros on $\T$.} of $c(\theta)$ on $\T$ by $c^{-1}(0)=\{\theta_1,\cdots,\theta_m\}$. For any $\beta$ with  $0<2\beta<\beta(\alpha)$ and $0<\delta<1$, there is there is a sequence $q_n\to \infty$ and a full Lebesgue measure set $\Theta=\Theta(\alpha,c^{-1}(0))$ such that 
for any $\theta\in\Theta$,   $c(\theta+n\alpha)$ satisfies:

\begin{align}\label{eq:Lambda-bound-Lip-analytic-0}
\min_{|k|\le q_{n}^{-1}e^{\delta\beta\,q_n}}\prod_{j=kq_n}^{(k+1)q_n-1}  |c(\theta+j\alpha)|>e^{-2\delta^2 \min\{\beta,1\}\,q_n},
\end{align}
\begin{align}\label{eq:Lambda-bound-Lip-analytic}
\min_{|m|\le e^{\delta\beta\,q_n}}\prod_{j=m}^{m+q_n-1}  |c(\theta+j\alpha)|>e^{-6\delta^2 \min\{\beta,1\}\,q_n}.
\end{align}

\end{cor}

\proof
 Clearly, there is an analytic function $\widetilde g(\theta)$ such that:
\begin{align}\label{eq:c-analytic}
c(\theta)=\widetilde g(\theta)\prod_{{\ell}=1}^m\big(e^{2\pi{\rm i}\theta}-e^{2\pi{\rm i}\theta_{\ell}}\big),\ \inf_\T|\widetilde g(\theta)|>0. 
\end{align}

Direct computation shows $\int_\T \ln|\widetilde g(\theta)|{\rm d}\theta=\int_\T \ln|c(\theta)|{\rm d}\theta=0 $ and 
\begin{align} \label{eq:c-0}
c(\theta)=\widetilde g(\theta)\prod_{{\ell}=1}^m\big(e^{2\pi{\rm i}\theta}-e^{2\pi{\rm i}\theta_{\ell}}\big)
=\widetilde g(\theta)(2{\rm i})^m\prod_{{\ell}=1}^me^{{\rm i}\pi(\theta+\theta_{\ell})}\,\sin\pi(\theta-\theta_{\ell}).
\end{align}
Therefore, 
\begin{align}
\frac{|c(\theta)|}{\prod_{{\ell}=1}^m|\sin\pi(\theta-\theta_{\ell})|}
=2^m|\widetilde g(\theta)|.
\end{align}

Apply Lemma \ref{lem:LipBound} to (\ref{eq:c-0}) where $\tau_1=\cdots=\tau_m=1$ and $|g(\theta)|\equiv2^m|\widetilde g(\theta)|$, we have (\ref{eq:Lambda-bound-Lip}) and (\ref{eq:cmqn}) hold with $\Lambda_1=2\delta^2\min\{\beta,1\}$ and $\Lambda=6\delta^2\min\{\beta,1\}$. 
\qed

\section{Spectral continuity: proof of Theorem \ref{thm:speconti}}\label{sec:spe-conti}
Following the notations and assumptions in Theorem \ref{thm:speconti}, consider \begin{equation}\label{eq:Jacobi}
(Hu)_n=w_n u_{n+1} +\overline{w}_{n-1} u_{n-1}+v_n u_n,\ n\in\Z.
\end{equation}
Assume that there are positive constants $\beta,\delta,\Lambda>0$ and a sequence of positive integers $q_n\to \infty$ such that $w_n,v_n$ has $\beta$-$q_n$ almost periodicity and $w_n$ has $(\Lambda,\beta)$-$q_n$ bound.

The key observation is: if $H$ has $\beta$-$q_n$ almost periodicity, then it can be approximated by a $q_n$ periodic operator exponentially fast in a finite (exponentially large) lattice. The estimates on the $q_n$ periodic operator eventually lead to the quantitative upper bound for the $m$-function as in  (\ref{def:speconti}) through the subordinacy theory Lemma \ref{lem:JLineq}. 

In view of Lemma \ref{lem:JLineq}, let $v^{\varphi}$ be the right half line solution to $Hu=Eu$ with initial condition $\varphi$ and $\ell=\ell(\varphi,\varepsilon,E)$ is defined as in (\ref{eq:subor-l}). As a direct consequence of Lemma \ref{lem:JLineq} and Proposition \ref{prop:DKL}, the following relation between the power law of $\|v^{\varphi}\|_{\ell}$ and the spectral continuity was proved in \cite{JZ} (see Lemma 2.1 and the proof of Theorem 6 there). 

\begin{lemma}\label{lem:Powerlaw}
	Fix $0<\gamma<1$. Suppose for $\mu$-a.e. $E$, there is a sequence of positive numbers $\eta_k\to 0$ and $L_k=\ell_k(\varphi,\eta_k,E)\to\infty$ such that
	for any $\varphi$ 
	\begin{equation}\label{eq:uvLk}
	1/16\big(L_k\big)^{\gamma} \le \|v^{\varphi}\|^2_{L_k}\le  \big(L_k\big)^{2-\gamma}.
	\end{equation}
	
	Then the spectral measure $\mu$ is $\gamma$-spectral continuous. 
\end{lemma}

Let  $A(n;E)$ be defined as in (\ref{def:An}). Denote by ${\rm Tr}\,A$ the trace of any matrix $A\in{\rm GL}(2,\mathbb{C})$. The following estimate on ${\rm Tr}\,A(q_n;E)$ is the key to prove the above power law and spectral continuity. 
\begin{thm}\label{thm:trace}
Let $H, \beta,\delta,\Lambda$ and $q_n$ be given as in (\ref{eq:Jacobi}). Suppose $\beta>260(1+\frac{1}{\delta})\Lambda$, then for  $\mu\ a.e.\ E$, there exists $K(E)\in\N$, for $k\geq K(E)$, we have
\begin{align} \label{eq:trace}
|Tr A(q_k;E)| <2-2e^{-60\Lambda q_k}.
\end{align}
For any $0<\gamma<1$, assume further that \begin{align}\label{eq:beta>Lambda}
\beta>300(1+\frac{1}{\delta})\frac{\Lambda}{1-\gamma},
\end{align}
we have the power law required by (\ref{eq:uvLk}). 
\end{thm}
Let 
\begin{align}
C=C(\delta,\Lambda)=300(1+\frac{1}{\delta})\Lambda.
\end{align}
Combining Lemma \ref{lem:Powerlaw} and  (\ref{eq:beta>Lambda}) in Theorem \ref{thm:trace}, 
if $\beta>C$, then $\mu$ is $\gamma$-spectral continuous for any $\gamma<1-\frac{C}{\beta}<1$ and therefore ${\rm dim}_{\rm spe}(\mu)\ge 1-\frac{C}{\beta}$. This proves Theorem \ref{thm:speconti}.\\

The trace estimate (\ref{eq:trace}) shows that spectrally almost everywhere, $A(q_k;E)$ is strictly elliptic eventually. The quantitative estimate (\ref{eq:trace}) allows us to iterate the transfer matrix up to the length scale $e^{\Lambda q_k}$ , which gives a well control on the norm of $A({q_k};E)$. The norm estimate eventually leads to the power law as required in (\ref{eq:uvLk}) through (\ref{eq:eigeneq}).

The proof of (\ref{eq:beta>Lambda}) and the required power law follows the outline of the Schr\"odinger case (see \cite{JZ}, Lemma 2.1). The main difference is now the transfer matrix $A(n;E)$ is in $GL(2,\mathbb{C})$. We need to consider some transformations introduced in \cite{HJ} which conjugate $A(n;E)$ to some ${\rm SL}(2,\R)$ matrix. Then many important techniques developed in \cite{JZ} for ${\rm SL}(2,\R)$ cocycles are now applicable. The trace estimate (\ref{eq:trace}) leads to a norm estimate of $A({q_k};E)$ and eventually leads to the estimate (\ref{eq:uvLk}) for the truncated ${\ell}^2$ norm of the eigenfunction $v^{\varphi}$ by (\ref{eq:eigeneq}). We will omit the details here and focus on the proof of the trace estimate (\ref{eq:trace}). For the sake of completeness, we sketch the proof of (\ref{eq:beta>Lambda}) and the power law (\ref{eq:uvLk}) in the Appendix \ref{A:pflemPower} for reader's convenience.

The rest of the section is organized as follows: In section \ref{sec:Conj-rT}, we introduce the transformation we will use to conjugate ${\rm GL}(2,\mathbb{C})$ to ${\rm SL}(2,\R)$ and develop all the useful lemmas about the conjugate. In section \ref{sec:hyp}, we study the case where the trace of the transfer matrix is greater than 2. In section \ref{sec:ellip}, we study the case where the trace of the transfer matrix is close to 2. 

Throughout this section, we assume $v_n,w_n$ have $\beta$-$q$ almost periodicity and $w_n$ has $(\Lambda,\beta)$-$q$ bound for some $q$ large enough such that $e^{-(\beta-2\Lambda)q}<1/10$.
We also use the induced estimates (\ref{eq:wrm-lower})-(\ref{eq:Arm-upper}) discussed in section \ref{sec:MoreLambdaBound} directly, refered also as $\beta$-$q$ almost periodicity and $(\Lambda,\beta)$-$q$ bound.

\subsection{Conjugate between ${\rm SL}(2,\R)$ and $GL(2,\mathbb{C})$ matrices.}\label{sec:Conj-rT}
The trace estimate (\ref{eq:trace}) was first proved in \cite{JZ} for ${\rm SL}(2,\R)$ cocycles. The generalization to $GL(2,\mathbb{C})$ case is very delicate. We need to consider the following transformation: let \begin{align} \label{def:T}
T_n= \left (\begin{matrix}  1   &  0 \\ 0  &   \sqrt{\frac{w_{n}}{\overline{w}_{n}}} \end{matrix} \right )
\end{align} and
\begin{align} \label{def:r}
r_n= \frac{w_{n+1}}{\sqrt{|w_{n+1}\,w_{n}|}}.
\end{align}

Let $A_n(E)$ be given as in (\ref{def:AnDn}). 
Define \begin{align} \label{def:tildeA}
\widetilde{A}_n(E):=r_{n-1} \, T^{-1}_{n}\, A_n\, T_{n-1}=\frac{1}{\sqrt{|w_nw_{n-1}|}} \left (\begin{matrix}  E-v_n   &  -|{w}_{n-1}| \\ |w_n|  &   0 \end{matrix} \right ).
\end{align}
The n-transfer matrix $\widetilde A(n,m;E)$ and $\widetilde A(n;E)$ for $\widetilde A_n$ will be defined in the same way as in (\ref{def:Anm}):
\begin{equation}\label{def:tildeAnm}
\widetilde A(n,m;E)=\prod_{j=m}^{n+m-1} \widetilde A_j(E), \  \  n\in\N^+, \ m\in\Z
\end{equation}
and 
\begin{equation}\label{def:tildeAn}
\widetilde A(n;E)=\widetilde A(n,1;E), n>0;\  \widetilde A(0;E)=Id;\ \widetilde A(n;E)=\widetilde A^{-1}(-n,n+1;E), n<0.
\end{equation}
We also denote the scalar product of $r_n$ in the same way as $w(n,m)$ in (\ref{def:cn}) for $n\in\N^+$,
\begin{equation}\label{def:rn}
r(n,m)=\prod_{j=m}^{n+m-1}r_j,\ \ m\in\Z.
\end{equation}

Direct computation shows
\begin{align}\label{eq:conjA-A}
\widetilde A(n,m;E)=r(n,m-1) \, T^{-1}_{n+m-1}\, A(n,m;E)\, T_{m-1}.
\end{align}

In view of (\ref{def:T}) and \eqref{def:tildeA}, it is easy to check  that $\|T_n\|=1$ and $\widetilde{A}_n,\widetilde A(n,m)\in {\rm SL}(2,\R)$ for any $n,m$. By (\ref{def:tildeA}) and (\ref{eq:conjA-A}), we are able to apply the techniques developed in \cite{JZ} for ${\rm SL}(2,\R)$ matrix and then swtich between the singular ${\rm GL}(2,\mathbb{C})$ case and the ${\rm SL}(2,\R)$ case.

The $\beta$-almost periodicity and the $(\Lambda,\beta)$ boundedness of $w_n$ imply the $\beta$-almost periodicity of $r$ and $T$ in the following sense:

\begin{lem}\label{lem:rT}
If $\beta>2\Lambda$, then for all $m\in\Z$ such that $|m|<e^{\delta\beta q},$
\begin{align}
\Big||r^{\pm}(q,m)|-1\Big| <e^{-(\beta-2\Lambda) q}  \label{eq:r}\\
\|T^{-1}_{m+q} \cdot T_{m}-I \|=\|T_{m}\cdot T^{-1}_{m+q}  -I \| <4e^{-(\beta-2\Lambda) q}. \label{eq:T}
\end{align}

Assume further that $N\in\N^+,Nq\le e^{\delta \beta q}$, then
\begin{align}
\Big||r^{\pm}({Nq,0})|-1\Big| <N e^{-(\beta-2\Lambda) q}   \label{eq:rexpbeta}\\
\|T_0\cdot T^{-1}_{\- Nq}  -I \|=\|T^{-1}_{ Nq} \cdot T_0-I \| <4Ne^{-(\beta-2\Lambda) q} . \label{eq:Texpbeta}
\end{align}

\end{lem}
Note that $r(n,m)$ and $T_n$ are essentially scalar products, the proof is based on the following direct computation:
\proof
 Set $z_m=\frac{w_{m}}{w_{m+q}}$. 
By (\ref{eq:wq/wq}), for $|m|\le e^{\delta \beta q}$ and $q$ large,
\begin{align*}
\big||z_m|^{\pm}-1\big|\le|z^{\pm}_m-1|<e^{-(\beta-2\Lambda) q}<\frac{1}{2} .
\end{align*}

Clearly, $|r_n|=\sqrt{\frac{|w_{n+1}|}{|w_{n}|}}$. In view of (\ref{def:rn}) and \eqref{def:T},  we have
\begin{align}
\Big||r(q,m)|-1\Big| = \Big|\sqrt{\frac{|c_{m+1}|}{|c_{m}|}  \cdots \frac{|c_{m+q}|}{|c_{m+q-1}|}}-1\Big| 
=\Big|\sqrt{|z_m|^{-1}}-1\Big| <e^{-(\beta-2\Lambda) q}  \label{eq:temp-rqm}
\end{align}

and
\begin{align*}
T^{-1}_{m+q}\,T_{m}=
 \left (\begin{matrix}  1   &  0 \\ 0  &   \sqrt{\overline{w}_{m+q}{w^{-1}_{m+q}}} \end{matrix} \right )
 \left (\begin{matrix}  1   &  0 \\ 0  &   \sqrt{w_{m}{\overline{w}^{-1}_{m}}} \end{matrix} \right ) 
 =\left (\begin{matrix}  1   &  0 \\ 0  &  \frac{w_{m}w_{m+q}^{-1}}{|w_{m}w_{m+q}^{-1}|} \end{matrix} \right )
 =\left (\begin{matrix}  1   &  0 \\ 0  & \frac{z^{}_m}{|z^{}_m|} \end{matrix} \right )
\end{align*} 
Therefore,
\begin{align*}
\|T^{-1}_{m+q}\,T_{m}-I\|\le \Big|\frac{z^{}_m}{|z^{}_m|}-1\Big|=\frac{\Big|z^{}_m-|z^{}_m|\Big|}{|z^{}_m|}\le 
\frac{2\big|z^{}_m-1\big|}{|z^{}_m|}\le 4e^{-(\beta-2\Lambda) q} .
\end{align*} 

In particular, in (\ref{eq:temp-rqm}), let $m=0,q,2q,\cdots,(N-1)q$ for $Nq<e^{\delta \beta q}$. Direct computation shows that 
\begin{align*}
\Big||r(Nq,0)|-1\Big| = \Big|\prod_{k=0}^{N-1}|r(q,kq)|-1\Big| 
<Ne^{-(\beta-2\Lambda) q}, 
\end{align*}

\begin{align*}
T^{-1}_{Nq}\,T_{0}=
 \left (\begin{matrix}  1   &  0 \\ 0  &   \sqrt{\overline{w}_{Nq}{w^{-1}_{Nq}}} \end{matrix} \right )
 \left (\begin{matrix}  1   &  0 \\ 0  &   \sqrt{w_{0}{\overline{w}^{-1}_{0}}} \end{matrix} \right ) 
 =\left (\begin{matrix}  1   &  0 \\ 0  &  \frac{w_{0}w_{Nq}^{-1}}{|w_{0}w_{Nq}^{-1}|} \end{matrix} \right ),
\end{align*} 
and 
\begin{align*}
\|T^{-1}_{Nq}\,T_{0}-I\|\le \Big|\frac{w_{0}w_{Nq}^{-1}}{|w_{0}w_{Nq}^{-1}|}-1\Big|
&\le \Big|\prod_{k=0}^{N-1}\frac{z_{kq}}{|z_{kq}|}-1\Big|\\
&\le \sum_{k=0}^{N-1}\Big|\frac{z_{kq}}{|z_{kq}|}-1\Big|\le4Ne^{-(\beta-2\Lambda) q}. 
\end{align*} 
\qed

(\ref{eq:conjA-A}) only implies $\|\widetilde A(q;E)\|\approx \| A(q;E)\| $, while $\|A(q;E)-\widetilde A(q;E)\|$ is not necessarily small. 
Lemma \ref{lem:rT} actually shows $\widetilde A(q;E)$ and $A(q;E)$ are close to each other up to a conjugate. This will be enough to control the difference between their traces. Fix $E$, we write $A(n,m)=A(n,m;E)$ for short. 
 \begin{lem}\label{lem:normA-A}
For all $m\in\Z$ such that $|m|<e^{\delta\beta q}$, let $\Phi={\rm Arg}\,r(q,m)$ be the Principal value of $r(q,m)\in\C$. For $\beta>4\Lambda$,
\begin{align}
\frac{1}{2}\|A(q,m)\|\le \|\widetilde{A}(q,m)\|、\le 2\|A(q,m)\| < 2e^{2\Lambda q} \label{eq:normA}
 \end{align}

 \begin{align}
 \|\widetilde{A}(q,m+1)-e^{{\rm i}\Phi}\,T^{-1}_{m}\, A(q,m+1)\, T_{m} \| < 12e^{-(\beta-4\Lambda)q} \label{eq:normA-A}
 \end{align}
 and consequently,
  \begin{align}
\Big| |Tr \widetilde{A}(q,m)| -|Tr A(q,m)| \Big|<  12e^{-(\beta-4\Lambda)q} \label{eq:tr-tr}.
 \end{align}
 \end{lem}

\proof
By (\ref{eq:r}), we have $|r^{\pm}(q,m)|\le 2$. (\ref{eq:normA}) follows from (\ref{eq:Arm-upper}) and (\ref{eq:conjA-A}) since $\|T^{\pm}_m\|=1$.
By (\ref{eq:conjA-A}), we have 
 \begin{align}
\widetilde{A}(q,m+1)=r(q,m)T^{-1}_{m+q} A(q,m+1) T_m =\big(|r(q,m)|T^{-1}_{m+q}T_m\big)\,e^{{\rm i}\Phi}\,T^{-1}_m A(q,m+1) T_m 
 \end{align}
Therefore,
 \begin{align}
\|\widetilde{A}(q,m+1)-e^{{\rm i}\Phi}\,T^{-1}_m A(q,m+1) T_m \|&=
\|\big(|r(q,m)|T^{-1}_{m+q}T_m-I\big)\,e^{{\rm i}\Phi}\,T^{-1}_m A(q,m+1) T_m  \|    \\
&\le \|\big(|r(q,m)|T^{-1}_{m+q}T_m-I\big)\| \cdot \|e^{{\rm i}\Phi}\,T^{-1}_m A(q,m+1) T_m  \|  \\
& \le 6e^{-(\beta-2\Lambda)q}\|A(q,m+1)\|  \\
&  \le 12e^{-(\beta-4\Lambda)q} . \label{eq:normA-AA}
 \end{align}
The last inequality follows from (\ref{eq:r}) and (\ref{eq:T}) since 
 \begin{align*}
\|\big(|r(q,m)|T^{-1}_{m+q}T_m-I\big)\| \le \big||r(q,m)|-1\big|\cdot \|T^{-1}_{m+q}T_m-I\|+\big||r(q,m)|-1\big|+\|T^{-1}_{m+q}T_m-I\|.
 \end{align*}
(\ref{eq:tr-tr}) follows directly from (\ref{eq:normA-A}) since $|{\rm Tr}\, A(q,m+1)|=\Big|{\rm Tr}\big(\,e^{{\rm i}\Phi}\,T^{-1}_m A(q,m+1) T_m\,\big)\Big|$.

\qed

Standard telescoping argument allows us to pass the $\beta$-almost periodicity from the sequences $w_n,v_n$ to the matrices $A(n,m),\widetilde A(n,m)$,  up to product length $q$,
\begin{lem}\label{lem:Abeta}
For all $m\in\Z$ such that $|m|<e^{\delta \beta q}$, $\beta>6\Lambda$, 
\begin{align}\label{eq:Abeta}
\|A(q,m;E)-A(q,m+q;E)\|\le e^{(-\beta+6\Lambda)q}
\end{align}
and 
\begin{align}\label{eq:tildeAbeta}
\|\widetilde A(q,m;E)-\widetilde A(q,m+q;E)\|\le e^{(-\beta+6\Lambda )q}.
\end{align}
\end{lem}

\proof
Write $m'=m+q$ for short. 
\begin{align*}
A(q,m)-A(q,m')&=\sum_{j=0}^{q-1}A(q-j-1,m+j+1) \Big(A_{m+j}-A_{m'+j}\Big)A(j,m') \\
&=\sum_{j=0}^{q-1}\frac{D(q-j-1,m+j+1)}{w(q-j-1,m+j+1)} 
\Big(\frac{D_{m+j}}{w_{m+j}}-\frac{D_{m'+j}}{w_{m'+j}}\Big)\frac{D(j,m')}{w(j,m')}.
\end{align*} 

By the trivial upper bound (\ref{def:D-upper}) for $D(n,m)$ and the lower bound (\ref{eq:wrm-lower}) for $w(n,m)$, we have 

\begin{align*}
\|A(q,m)-A(q,m')\|
&\le\sum_{j=0}^{q-1}\frac{e^{(q-j-1)\Lambda}}{|w(q-j-1,m+j)|} 
\Big|{w_{m'+j}}{D_{m+j}}-{w_{m+j}}{D_{m'+j}}\Big|\frac{e^{\Lambda j}}{|w(j+1,m')|}\\
&\le\sum_{j=0}^{q-1}\frac{e^{(q-j-1)\Lambda}}{e^{-2\Lambda q}} 
\Big|{w_{m'+j}}{D_{m+j}}-{w_{m+j}}{D_{m'+j}}\Big|\frac{e^{j\Lambda}}{e^{-2\Lambda q}}\\
&\le q\,e^{5\Lambda q}\max_{|m|\le e^{\delta \beta q}}\Big|{w_{m+q}}{D_{m}}-{w_{m}}{D_{m+q}}\Big|\\
&\le q\,e^{5\Lambda q}\max_{|m|\le e^{\delta \beta q}}
\Big(\big|({w_{m+q}}-w_m){D_{m}}\big|+\big|{w_{m}}({D_{m+q}}-D_m)|\Big) \\
& \le 2q\,e^{5\Lambda q}\, e^{\Lambda}e^{-\beta q} \\
& \le e^{-(\beta-6\Lambda) q}
\end{align*} 
provided $\sup_n|w_n|,\sup_{n,E}\|D_n\|\le e^{\Lambda}$ and $q$ large such that $2qe^{\Lambda}\le e^{\Lambda q}$.\\

 Let 
$\widetilde w_n=\sqrt{|w_nw_{n-1}|}$, $ \widetilde D_n=\left (\begin{matrix}  E-v_n   &  -|{w}_{n-1}| \\ |w_n|  &   0 \end{matrix} \right ).
$ Define $\widetilde w(n,m),\widetilde D(n,m;E)$ exact in the same as for $w,D$. It is easy to check that the $\Lambda$ bounds of $w_n$ and $D_n$ hold true for $\widetilde w_n,\widetilde D_n$:
\begin{align}
 |\widetilde w(n,m)| \leq e^{\Lambda n}, \ \sup_{E} \|\widetilde D(n,m;E)\| \leq e^{\Lambda n}, \ \ \forall n\ge0 ,m\in\Z
\end{align}
and 
 \begin{align}\label{eq:tildecboundlower}
 |\tilde w(r,m)| \ge e^{-2\Lambda\,q}, \ \ 0\le r\le q,\ |m|\leq e^{\delta\beta q}.
\end{align}

 The $\beta$-almost periodicity of $w_n,v_n$ are also passed directly to $\widetilde w_n,\widetilde D_n$:
\begin{align}
\max_{|m|\leq e^{\delta\beta q}} |\widetilde w_{m} -\tilde w_{m\pm q}|\leq 2e^{\Lambda}{e^{-\beta q}} ,
\quad \max_{|m|\leq e^{\delta\beta q}} \|\widetilde D_{m} -\widetilde D_{m\pm q}\|
\leq {e^{-\beta q}} .
\end{align}

By the definition of $\widetilde A_n$ in (\ref{def:tildeA}), we have 
$\widetilde A_n=\frac{1}{\widetilde w_n} \widetilde D_n$. Exact the same computation proves that 
\begin{align}
\|\widetilde A(q,m)-\widetilde A(q,m')\|\le q(2e^{2\Lambda}+e^{\Lambda})e^{5\Lambda q}\, e^{-\beta q}  \le e^{-(\beta-6\Lambda) q}.
\end{align}

\qed

The above  telescoping argument can not be extended to exponential scale  $e^{\Lambda q}$ as for $r$ and $T$ in (\ref{eq:rexpbeta}),(\ref{eq:Texpbeta}) directly. One main reason is we lose control of the matrix norm super-exponentially as $\|D(e^{\Lambda q})\|\lesssim e^{\Lambda e^{\Lambda q}}$. Such gowth can not be controlled by condition such as $\beta\gtrsim \Lambda$.  The key to prove the trace estimate (\ref{eq:trace}) is to avoid using such rough bound for the matrix norm at an exponential scale.  This is one breakthrough in \cite{JZ}. By all the above estimates of the conjugate $r,T$ and some simple linear algebra facts of ${\rm SL}(2,\R)$ matrix found in \cite{JZ}, we are able to prove this extension for the ${\rm GL}(2,\mathbb{C})$ case. We will see more details in the next two subsections. 

Similar to \cite{JZ}, we consider the following two cases where $|Tr A(q)|$ is away from 2 and close to 2.

\subsection{The case where the trace is away from 2.}\label{sec:hyp}
We start with the hyperbolic case in the following sense:
let 
\begin{align}\label{def: S1}
S^{1}_q=\{E: \ |{\rm Tr} A(q;E)| > 2+2e^{-60\Lambda q}\}
\end{align}
We may fix $E$ and write $A(q)=A(q;E)$ for simplicity whenever it is clear. 

\begin{lemma}\label{lem:S1}
Let $q_n$ be given as in Theorem \ref{thm:trace}. If $\beta>(260+\frac{61}{\delta})\Lambda$, then the set
\begin{align}\label{def: limitsup}
\limsup_{n\to\infty} \, S^{1}_{q_n}
=\Big\{E:\ E\ {\textrm{belongs to infinitely many}}\ S^{1}_{q_n}\Big\}
\end{align}
 has spectral measure zero.
\end{lemma}

Lemma \ref{lem:normA-A} implies for large $\beta$, ${\rm Tr} \widetilde{A}(q;E)$ and ${\rm Tr} {A}(q;E)$  lie in the same region, i.e., if $E\in S^{1}_q$, then
\begin{align}\label{eq:trAAqhyp}
|{\rm Tr} \widetilde{A}(q;E)|> 2+2e^{-60\Lambda q}-12e^{-(\beta-4\Lambda)q}>2+e^{-60\Lambda q}, 
\end{align}
provided $e^{(\beta-64\Lambda)q}>12$ . 

The following linear algebra facts were proved in \cite{JZ}
\begin{lemma}\label{lem:JZ2.7} Suppose $G\in {\rm SL}(2,\R)$ with
  $2<|{\rm Tr\,}G|\le 6.$ The invertible matrix $B$ such that
\begin{equation}\label{G0}
    G=B\left(
                \begin{array}{cc}
                  \rho & 0 \\
                  0 & \rho^{-1} \\
                \end{array}
              \right)B^{-1}
\end{equation}
where $\rho^{\pm1}$ are the two conjugate real eigenvalues of $G$ with $|{\rm det}B|=1$ satisfies
\begin{equation}\label{B0}
    \|B\|=\|B^{-1}\|<\frac{\sqrt{\|G\|}}{\sqrt{|{\rm Tr\,}G|-2}}
\end{equation}
If $|{\rm Tr\,}G|> 6$,  then $\|B\|\le\frac{2\sqrt{\|G\|}}{\sqrt{|{\rm Tr\,}G|-2}}$.
\end{lemma}
Apply the above lemma to $ \widetilde{A}(q;E) \in {\rm SL}(2,\R)$ satisfying (\ref{eq:normA}) and (\ref{eq:trAAqhyp}),we have the following decomposition 
\begin{equation}\label{eq:decomp}
    \widetilde{A}(q)=B\left(
                \begin{array}{cc}
                  \rho & 0 \\
                  0 & \rho^{-1} \\
                \end{array}
              \right)B^{-1}
\end{equation}
where $\rho^{\pm1}$ are the two conjugate real eigenvalues of $\widetilde{A}(q)$ with $|\rho|>|{\rm Tr}\widetilde A(q)|-1>1+e^{-60\Lambda q}$ and $B$ satisfies $|{\rm det}B|=1$ and
\begin{equation}\label{eq:decompB}
    \|B\|=\|B^{-1}\|<e^{32\Lambda q}.
\end{equation}

By (\ref{eq:decomp}) and (\ref{eq:decompB}), we have that for any $N\in\N^+$,
\begin{align} \label{eq:tildeANqbound}
\widetilde A^N(q):=[\widetilde A(q)]^N=B\left(
                \begin{array}{cc}
                  \rho^N & 0 \\
                  0 & \rho^{-N} \\
                \end{array}
              \right)B^{-1},\ \    \|\widetilde A^N(q)\|\le e^{64\Lambda q}|\rho|^N
\end{align}

In the rest of this section, consider $\beta>\frac{61}{\delta}\Lambda$ and set 
\begin{align}\label{def:N-hyp}
N=[e^{61\Lambda q}]<e^{\delta \beta q}. 
\end{align}
The above decomposition now turns the matrix product $[\widetilde A(q)]^N$ into a scalar product of $\rho^N$ with a uniformly controlled conjugate $B$ (independent of $N$).  This is one key algebra ingredient observed in \cite{JZ}. This technique now allows us to extend the orbit of $\widetilde{A}(q)$ to the exponentially long scale $N=e^{61\Lambda q}$.

The following technique lemma was proved in \cite{JZ} (see Lemma A.1 there):
\begin{lemma}\label{lem:JZA1}
Suppose $G$ is a two by two matrix satisfying
\begin{equation}\label{G^j}
    \|G^j\|\le M<\infty, \ \ \ \textrm{for all}\ 0<j\le N\in\N^+,
\end{equation}
where $M\ge 1$ only depends on $N$. Let $G_j=G+\Delta_j$,
$j=1,\cdots,N,$ be a sequence of two by two matrices with
\begin{equation}\label{delta}
    \delta=\max_{1\le j\le N}\ \|\Delta_j\|.
\end{equation}
If
\begin{equation}\label{deltasmall}
    NM\delta<1/2,
\end{equation}
then for any $1\le n\le N$
\begin{equation}\label{G-product}
    \|\coprod_{j=1}^nG_j-G^n\|\le 2NM^2\delta.
\end{equation}
\end{lemma}

Let $N=[e^{61\Lambda q}]$, $G=\frac{1}{\rho}\widetilde A(q)$ and $G_j=\frac{1}{\rho}\widetilde A(q,jq+1)$, $|j|=0,1,\cdots N$.  By (\ref{eq:normA})
and (\ref{eq:tildeAbeta}), it is easy to check that $\|G^j\|\le e^{64\Lambda q} $ and $\|G_j-G\|\le Ne^{(-\beta+6\Lambda) q}\le e^{(-\beta+67\Lambda)q}$. The above lemma is applicable provided $\beta>(260+\frac{61}{\delta})\Lambda$.  One can prove that
\begin{align}
\label{eq:NqNq}
    &\|\widetilde{A}({Nq})-\widetilde{A}^N(q)\| \leq |\rho|^N\,e^{(-\beta+260\Lambda)q}
    \\
\label{eq:Nq-Nq}
   & \|\widetilde{A}(-{Nq})-\widetilde{A}^{-1}({Nq})\| \leq 2|\rho|^N\,e^{(-\beta+260\Lambda)q}
\end{align}
The proof (\ref{eq:NqNq}) and (\ref{eq:Nq-Nq}) is a direct application of Lemma \ref{lem:JZA1} and resembles the proof of Claim 3, \cite{JZ}. We omit the details here.

Similar to (\ref{eq:normA-A}), we can prove $A(\pm Nq)$ and $\widetilde A(\pm Nq)$ are close to each other up the size $|\rho|^N$. 
\begin{lem}\label{lem:normANq-ANq}
Let $\eta=r^{-1}(Nq,0), \zeta=r(Nq,-Nq)$ and $\phi={\rm Arg}\,\eta, \psi={\rm Arg}\,\zeta$ be the Principal values of $\eta$ and $\zeta$ accordingly. 
For $\beta>(260+\frac{61}{\delta})\Lambda$,

 \begin{align}\label{eq:A+Nq}
 \|{A}^{\pm}(Nq)-e^{\pm{\rm i}\phi}\,T_{0}\, \widetilde A^{\pm}(Nq)\, T^{-1}_{0} \| <  e^{(-\beta+127\Lambda)q}|\rho|^N,
 \end{align}

 \begin{align}\label{eq:A-Nq}
 \|{A}(-Nq)-e^{{\rm i}\psi}\,T_{0}\, \widetilde A(-Nq)\, T^{-1}_{0} \| <  e^{(-\beta+127\Lambda)q}|\rho|^N,
 \end{align}
and consequently,
 \begin{align}\label{eq:A-Nq-Nq}
 \|{A}^{-1}(Nq)-e^{-{\rm i}(\phi+\psi)}\,  A(-Nq) \| <  4e^{(-\beta+260\Lambda)q}|\rho|^N.
 \end{align}

\end{lem}

\proof
By (\ref{def:tildeA}), 
\begin{align}
{A}(Nq)=\eta\, T_{Nq}\, \widetilde A(Nq)\, T^{-1}_0 =\big(|\eta|\,T_{Nq}T^{-1}_0\big)\,e^{{\rm i}\phi}\,T_{0}\, \widetilde A(Nq)\, T^{-1}_{0}.
 \end{align}
Therefore, by (\ref{eq:rexpbeta}) and (\ref{eq:Texpbeta}), 
 \begin{align*}
\|{A}(Nq)-e^{{\rm i}\phi}\,T_{0}\, \widetilde A(Nq)\, T^{-1}_{0} \|&\le 
\big\|\,|\eta|\,T_{Nq}T^{-1}_0-I\big\|\cdot \|e^{{\rm i}\phi}\,T_{0}\, \widetilde A(Nq)\, T^{-1}_{0}  \|     \\
&\le 6e^{(-\beta+63\Lambda)q} \, \| \widetilde A(Nq) \|  \\
& \le  e^{(-\beta+127\Lambda)q}|\rho|^N,
 \end{align*}
 provided $e^{(\beta-260\Lambda)q}>1$ and $e^{\Lambda q}>12$. The last inequality follows from (\ref{eq:tildeANqbound}) and (\ref{eq:NqNq}):
 \begin{align*}
 \| \widetilde A(Nq) \| \le 
\|\widetilde{A}^N(q)\|+ e^{(-\beta+260\Lambda)q}\,  |\rho|^N 
\le e^{64\Lambda q}|\rho|^N+ e^{(-\beta+260\Lambda)q}\,|\rho|^N
\le 2e^{64\Lambda q}\, |\rho|^N.
 \end{align*}

 Note that $ \widetilde A(Nq)\in {\rm SL}(2,\R)$, then $ \| \widetilde A^{-1}(Nq) \|= \| \widetilde A(Nq) \|\le 2e^{64\Lambda q}|\rho|^N$. 
The proof for ${A}^{-1}(Nq)$ is exactly the same since 
\begin{align}
{A}^{-1}(Nq)=\eta^{-1}\, T_0\, \widetilde A^{-1}(Nq)\, T^{-1}_{Nq} 
=e^{-{\rm i}\phi}\,T_{0}\, \widetilde A^{-1}(Nq)\, T^{-1}_{0}\,\big(|\eta|^{-1}\,T_0T^{-1}_{Nq}\big).
 \end{align}
 
 (\ref{def:tildeA}) and (\ref{def:tildeAn}) imply that 
  \begin{align*}
{A}(-Nq)=A^{-1}(Nq,-Nq+1)&=\Big[r^{-1}(Nq,-Nq)T_0\widetilde A(Nq,-Nq+1)T^{-1}_{-Nq}\Big]^{-1}\\
&=r(Nq,-Nq)T_{-Nq}\widetilde A(-Nq)T^{-1}_0   .
 \end{align*}
 
   (\ref{eq:Nq-Nq}) implies that 
   $ \| \widetilde A(-Nq) \|\le \| \widetilde A^{-1}(Nq) \|+2|\rho|^N\,e^{(-\beta+260\Lambda)q}\le 3e^{64\Lambda q}|\rho|^N$.  Now by (\ref{eq:rexpbeta}) and (\ref{eq:Texpbeta}), exact the same argument for (\ref{eq:A+Nq}) proves (\ref{eq:A-Nq}) provided $e^{(\beta-260\Lambda)q}>2$ and $e^{\Lambda q}>18$. \\
 
The proof of (\ref{eq:A-Nq-Nq}) follows directly from (\ref{eq:Nq-Nq}), (\ref{eq:A+Nq}) and  (\ref{eq:A-Nq}) since
  \begin{align*}
 \|{A}^{-1}(Nq)-e^{-{\rm i}(\phi+\psi)}\,  A(-Nq) \|\le&  \|{A}^{-1}(Nq)-e^{-{\rm i}\phi}\,T_{0}\, \widetilde A^{-1}(Nq)\, T^{-1}_{0} \|          \\
&+\|e^{-{\rm i}\phi}\,T_{0}\, \widetilde A^{-1}(Nq)\, T^{-1}_{0}-e^{-{\rm i}\phi}\,T_{0}\, \widetilde A(-Nq)\, T^{-1}_{0}\|                \\
&+\|e^{-{\rm i}\phi}\,T_{0}\, \widetilde A(-Nq)\, T^{-1}_{0}-e^{-{\rm i}(\phi+\psi)}\,  A(-Nq)\| \\
\le & \|{A}^{-1}(Nq)-e^{-{\rm i}\phi}\,T_{0}\, \widetilde A^{-1}(Nq)\, T^{-1}_{0} \|  \\
&+\|\widetilde A^{-1}(Nq)-\widetilde A(-Nq)\|                \\
&+\|e^{{\rm i}\psi}\,T_{0}\, \widetilde A(-Nq)\, T^{-1}_{0}- A(-Nq)\|  \\
\le & 2 e^{(-\beta+127\Lambda)q}\,|\rho|^N+2 e^{(-\beta+260\Lambda)q}\,|\rho|^N \\
\le & 4e^{(-\beta+260\Lambda)q}\,|\rho|^N.
 \end{align*}
\qed

With the above preparation, we are in the place to prove Lemma \ref{lem:S1}. 
It is easy to see that all the estimates from (\ref{eq:NqNq}) to (\ref{eq:A-Nq-Nq}) preserve errors between the traces. 
Now combine (\ref{eq:NqNq}) with (\ref{eq:A+Nq}), we have 
\begin{align}
\Big||{\rm Tr}\,A(Nq)|-|{\rm Tr}\,\widetilde A^N(q)|\Big|\le 2e^{(-\beta+260\Lambda)q}\,|\rho|^N\le \frac{1}{2}|\rho|^N,
\end{align}
provided $e^{(\beta-260\Lambda)q}>4$.  Therefore, by (\ref{eq:tildeANqbound}), 
\begin{align}\label{eq:traceANqlower}
|{\rm Tr}\,A(Nq)|\ge |{\rm Tr}\,\widetilde A^N(q)|- \frac{1}{2}|\rho|^N
\ge  \frac{1}{2}|\rho|^N.
\end{align}
(\ref{eq:A-Nq-Nq}) implies that for any vector $X\in\mathbb{C}^2$,
\begin{align}\label{eq:normA-Nq}
\|A^{-1}(Nq)X\|\le \|A(-Nq)X\| +4e^{(-\beta+260\Lambda)q}\,|\rho|^N\|X\|\le  \|A(-Nq)X\| +\frac{1}{8}|\rho|^N\|X\|,
\end{align}
provided $e^{(\beta-260\Lambda)q}>32$.

By (\ref{def:tildeA}) and (\ref{eq:conjA-A}), it is easy to check that $|{\rm det}A(Nq)|=|r^{-1}(Nq,0)|$. 
Therefore, (\ref{eq:rexpbeta}) implies that
\begin{align}\label{eq:detANq}
|{\rm det}A(Nq)|<1+e^{-(\beta-63\Lambda) q}<2
\end{align}

Consider  the generalized eigenequation $Hu=Eu$ with normalized initial value
$X=\left(
  \begin{array}{c}
    u_1 \\
    u_0 \\
  \end{array}
\right), \ \ \|X\|=1$. 
By (\ref{eq:eigeneq}) and the Cayley-Hamilton theorem for ${\rm GL}(2,\mathbb{C})$ matrix $A(Nq)$, we have:
\begin{equation}\label{eq:uq}
   A({N q})\, X=\left(
  \begin{array}{c}
    u_{Nq+1} \\
    u_{Nq} \\
  \end{array}
\right), \ \ \ A({-Nq})\, X=\left(
  \begin{array}{c}
    u_{-Nq+1} \\
    u_{-Nq} \\
  \end{array}
\right).
\end{equation}
and 
\begin{align} \label{eq:CalyeyH}
A({N q})\, X+({\rm det}A)\cdot A^{-1}({N q})\, X=-({\rm Tr}\,A(Nq))\, X
\end{align}
Combine (\ref{eq:traceANqlower}),(\ref{eq:normA-Nq}),(\ref{eq:detANq}) with (\ref{eq:CalyeyH}), we have
\begin{align} \label{eq:normAX}
\|A({N q})\, X\|+\| A^{}(-{N q})\, X\|\ge \frac{1}{8}|\rho|^N. 
\end{align}

Now by the choice of $\rho$ and $N$, for $q$ large, we have
\begin{align} \label{eq:normAX-2}
\|A({N q})\, X\|+\| A(-{N q})\, X\|\ge \frac{1}{8}(1+e^{-60\Lambda q})^{[e^{61\Lambda q}]}\ge 4e^q
\end{align}
which implies
\begin{align} \label{eq:4uNq}
\max\big\{ |u_{Nq+1}|,\,  |u_{Nq}| ,\, |u_{-Nq+1}| ,\,  |u_{-Nq}| \big\}\ge e^q.
\end{align}

In conclusion, we can claim the existence of a subsequence of $u_n$ at energy $E$ with following exponential growth:
\begin{claim}\label{clm:Gordon}
Assume $v_n,w_n$ have $\beta$-$q$ almost periodicity as in (\ref{eq:a-beta}) and $w_n$ has $(\Lambda,\beta)$-$q$ bound (\ref{def:Lambda-bound}),(\ref{eq:Lambda-bound}) for $q>q_0(\Lambda,\delta,\beta)$. Suppose 
$E\in S^1_q$ and $\beta> (260+\frac{61}{\delta})\Lambda$, then there are integer sequences $x_q^1, x_q^2, x_q^3, x_q^4 \in \Z$ independent of $E$, such that $\min_i|x_q^i| \rightarrow \infty$ as $q \rightarrow \infty$ and
\begin{align}
\max_{i} |u_{x_q^i}^E | > e^q,
\end{align} 
where $u_n^E$ solves the half-line problem $Hu=Eu$ with normalized boundary condition $|u_0|^2+|u_1|^2=1$.
\end{claim}

Now Lemma \ref{lem:S1} follows directly from Claim (\ref{clm:Gordon}) and the following lemma:
\begin{thm}[Extended Schnol's Theorem, Lemma 2.4, \cite{JZ}]
Fix any $y>1/2$. For any sequence $|x_k|\to \infty$(where the sequence
is independent of $E$), for spectrally a.e. $E$, there is a generalized eigenvector $u^E$ of $Hu=Eu$, such that
$$|u^E_{x_k}|<C(1+|k|)^{y}. $$
\end{thm}

\subsection{The case where the trace is close to 2.}\label{sec:ellip}
In this part, we consider those energy $E$ where the trace of $A(q;E)$ is close to 2. 
Let 
\begin{align}
S^2_q=\{E: \big||{\rm Tr} A(q;E)| - 2\big| < 2e^{-60\Lambda q} \} \label{def:S2}
\end{align} 
Again we assume that $q$ is large and $v_n,w_n$ satisfy $\beta$-$q$ almost periodicity (\ref{eq:a-beta}) and $\Lambda$-$q$ bound in (\ref{eq:Lambda-bound}) with positive finite parameters $\beta,\Lambda,\delta$. We can prove that 
\begin{lemma}\label{lem:S2}
If $\beta> (130+\frac{29}{\delta})\Lambda$, then
\begin{equation}\label{eq:Sqmes}
    \mu(S_q)<e^{-\frac{1}{10}\Lambda q},
\end{equation}
where $\mu$ is the spectral measure of $H$.
\end{lemma}

{\noindent \bf Proof of Theorem \ref{thm:trace}:} Assume now $\beta>260(1+\frac{1}{\delta})\Lambda$.  Let $q_n$ be given as in Theorem \ref{thm:trace}. 
Lemma \ref{lem:S1} implies that for spectrally a.e. $E$, there is $K_1(E)$ such that,
\begin{equation}\label{eq:K1}
    |{\rm Tr\,}\, A(q_k;E)|<2+2e^{- 60\Lambda q_k}, \ \ \forall k\ge K_1(E)
\end{equation}

Combine Lemma \ref{lem:S2} with the Borel-Cantelli lemma, we have 
$
\mu\Big(\limsup_nS^2_{q_n}\Big)=0
$,
i.e.,  for spectrally a.e. $E$, there is $K_2(E)$ such that
\begin{equation}\label{eq:K2}
    \Big||{\rm Tr\,}\, A(q_k;E)|-2\Big|> 2e^{- 60\Lambda q_k}, \ \ \forall k\ge K_2(E). 
\end{equation}

Clearly, (\ref{eq:K1}) and (\ref{eq:K2}) complete the proof of Theorem \ref{thm:trace} by taking $K=\max\{K_1,K_2\}$. \qed \\

In the rest of the section, we focus on proving (\ref{eq:Sqmes}).  Similar to the hyperbolic case, ${\rm Tr}\, \widetilde{A}(q;E)$ and ${\rm Tr}\,A(q;E)$  are close up to exponential error by Lemma \ref{lem:normA-A}. More precisely, let 
\begin{align}\label{def:tildeS2}
\widetilde S^2_q:=\Big\{E:\ \big||{\rm Tr} \widetilde A(q;E)| - 2\big| < 3e^{-60\Lambda q}\, \Big\}
\end{align}
Clearly, Lemma \ref{lem:normA-A} implies that for $\beta>6\Lambda$, $S^2_q\subset \widetilde S^2_q$. \\

The following elementary  linear algebra facts were proved in \cite{JZ}
\begin{lemma}[Lemma 2.9, Lemma 2.10 \cite{JZ}]\label{lem:JZ2.9}
Suppose $A\in {\rm SL}(2,\R)$ has eigenvalues $\rho^{\pm 1}$, $\rho>1$. For any $k\in\N$, if ${\rm Tr\,}A\neq2$, then
\begin{equation}\label{eq:JZ2.9-1}
    A^k=\frac{\rho^k-\rho^{-k}}{\rho-\rho^{-1}}\cdot\Big(
       A-\frac{{\rm Tr}A}{2}\cdot I\Big)+\frac{\rho^k+\rho^{-k}}{2}\cdot I
\end{equation}
Otherwise,
$A^k=k(A-I)+I$.

Assume further that $\big||{\rm Tr}\,A|-2\big|<\tau<1$, then there are universal constants
  $1<C_1<\infty,c_1>1/3$ such that for  $1\le k\le \tau^{-1}$, we have
\begin{equation}\label{eq:JZ2.9-2}
    c_1<\frac{\rho^k+\rho^{-k}}{2}<C_1\ ,\ \   c_1k<\frac{\rho^k-\rho^{-k}}{\rho-\rho^{-1}}<C_1k. 
\end{equation}
\end{lemma}

Now fix $E\in\widetilde S^2_q$, the above lemma actually shows that the $k$-th power of $\widetilde{A}({q;E})$ grows almost linearly with respect to $k$ as :
\begin{align}\label{eq:Alinear}
\widetilde{A}^k (q)\sim k \big(\widetilde{A}(q)-\frac{1}{2}{\rm Tr}\, \widetilde{A}(q)\big)+I, \ \ 1\le k \le N.  
\end{align}
 This simple observation will be an important part of our quantitative estimates in the near parabolic case. The arguments to derive (\ref{eq:Sqmes}) from (\ref{eq:Alinear}) follow the outline of the near parabolic case in \cite{JZ} with slight modification concerning all the estimates of the conjugacy  in section \ref{sec:Conj-rT}. We sketch the proof below for reader's convenience.   \\ 

\noindent{\bf Proof of Lemma \ref{lem:S2}}: 
First, Lemma \ref{lem:JZ2.9}  provides the following norm estimates: 
there is absolute constant $C_1>0$ such that for any $1\le j<N=[e^{29 \Lambda q}]<e^{29\Lambda q}<e^{\delta \beta q}$, 
\begin{align}
  \|\widetilde A^j(q;E)\|<3C_1\,j\,\cdot\|\widetilde A(q)\| .
\end{align}
By (\ref{eq:normA}) and the choice of $N$, we have 
\begin{align}\label{eq:normAlinear}
  \|\widetilde A^j(q;E)\|<6C_1\,j\,e^{2\Lambda q }<j\,e^{3\Lambda q } <e^{32\Lambda q}.
\end{align} 
In the same way as the proof of (\ref{eq:NqNq}) and (\ref{eq:Nq-Nq}), for any $1\le k \le N$, combine (\ref{eq:normAlinear}) with (\ref{eq:tildeAbeta}), we can apply Lemma \ref{lem:JZA1}  to obtain 
 \begin{align}\label{eq:A(kq)-Ak(q)}
 \|\widetilde{A}(kq) -\widetilde{A}^{k}(q) \| \leq  e^{(-\beta+130\Lambda) q}<1,
 \end{align}
provided $\beta>(130+\frac{29}{\delta})\Lambda$ and $q$ large. 

In view of (\ref{eq:JZ2.9-1}), (\ref{eq:Alinear}) and (\ref{eq:A(kq)-Ak(q)}), it is clear that $\widetilde{A}(kq)$ has the same linear expansion as in  (\ref{eq:Alinear}). Combine (\ref{eq:JZ2.9-1}), (\ref{eq:A(kq)-Ak(q)}) with  the conjugate relation:
\begin{align}\label{eq:AX-AX}
A({kq})\, X=r^{-1}({kq},0)\, T(kq) \widetilde{A}(kq) T^{-1}(0)\, X, \ \ X\in\C^2,
\end{align}
we can prove that:
\begin{claim} \label{clm:JZclm5}
For any $\varepsilon>0,E$ and $\varphi\in[0,2\pi)$, let $\ell=\ell(\varphi,\varepsilon,E),u^{\varphi},v^{\varphi}$ be given as in (\ref{eq:subor-l}). Suppose $E \in S^2_q,\varepsilon<e^{-29\Lambda q}$ and $\beta> (130+\frac{29}{\delta})\Lambda$, then 
\begin{align}\label{eq:ulqlower}
\|u^{\varphi}\|_{\ell}^2 > e^{\frac{1}{6}\Lambda q}.
\end{align}
\end{claim}

The proof of Claim \ref{clm:JZclm5} follows the outline of the proof of Claim 5 in \cite{JZ}. The key is to use the linear expression (\ref{eq:Alinear}) to control both the upper and lower bound of bound of $\|A(n)\,X\|$. The main difference is we need to consider the conjugacy (\ref{eq:AX-AX}) and switch between the orbits of $A(n)\,X$ and $\widetilde A(n)\,\tilde X$.  We omit the details here. For sake of completeness, we include the proof in Appendix \ref{A:pfclm5}.

We proceed to prove Lemma \ref{lem:S2} by Claim \ref{clm:JZclm5}.  ${\rm Tr}\, A(q;E)$ is a polynomial in $E$ with degree $q$. $S^2_q$ can be written as a union of at most $q$ band: $S^2_q=\bigcup_{j=1}^{q}{I}_j$. 
Note ${\rm Tr}\, \widetilde A(q;E)$ is also a polynomial in $E$ with degree $q$ with real coefficients, by Proposition \ref{A:JM}, we have $|\widetilde S^2_q|\le  C_2\sqrt {6e^{- 60\Lambda q}}$, where $C_2$ only depends on $\|w\|_\infty,\|v\|_\infty$. Then this gives us a uniform control on the width of each band $I_j$:
 \begin{equation}\label{eq:Ij}
       S^2_q=\bigcup_{j=1}^{q}{I}_j,\ \ \ \varepsilon^j _q:=|{I}_j|\le |S^2_q|\le |\widetilde S^2_q|\le  e^{- 29\Lambda q}.
\end{equation}

Now pick $E_j\in I_j\bigcap \sigma(H)\neq\emptyset$ to be the center in the sense that $I_j\subset(E_j-\varepsilon^j _q,E_j+\varepsilon^j _q)$. 
For any $\varphi$, let $u^{\varphi}(E_j)$ be the
right half line solution associated with the energy $E_j$. By Claim \ref{clm:JZclm5}, we have
\begin{align}\label{eq:uEj}
\|u^{\varphi}(E_j)\|_{\ell_q(j)}^2\ge e^{\frac{1}{6}\Lambda q}, \ \ j=1,\cdots,q
\end{align}
where $\ell_q(j)=\ell(\varphi,E_j,\varepsilon^j _q)$ is given as in (\ref{eq:subor-l}).

A direct consequence of (\ref{eq:uEj}) and the subordinacy theory Lemma \ref{lem:JLineq} is
\begin{align}\label{eq:mEj}
\varepsilon^j _q\cdot |m_{\varphi}(E_j+ i\varepsilon^j _q)|< \frac{5+\sqrt{24}}{2}\cdot e^{-\frac{1}{6}\Lambda q}, \ \ j=1,\cdots,q
\end{align}
Then by (\ref{def:Mfunction}) and (\ref{eq:DKL}), we have
\begin{align}\label{eq:muIj}
\mu(I_j)\le \sup_\varphi\, 2\varepsilon^j _q |m_{\varphi}(E_j+ i\varepsilon^j _q)|< (5+\sqrt{24}) e^{-\frac{1}{6}\Lambda q}, \ \ j=1,\cdots,q. 
\end{align}

Clearly, (\ref{eq:muIj}) completes the proof of Lemma \ref{lem:S2} provided $q(5+\sqrt{24}) e^{-\frac{1}{6}\Lambda q}\le e^{-\frac{1}{10}\Lambda q}$. \qed

\section{Spectral Singularity for analytic  quasiperiodic Jacobi operator}\label{sec:speSingular}
In this section, we focus on analytic quasiperiodic potential given by 
$v_n=v(\theta+n\alpha),w_n=c(\theta+n\alpha),n\in\Z,\theta\in\T$ where $v\in C^{\omega}(\T,\R)$ and $c\in C^{\omega}(\T,\mathbb{C})$ are analytic functions on $\T$ taken values in $\R$ and $\mathbb{C}$ respectively. Both $v(\theta)$ and $c(\theta)$ have bounded analytic extensions to the strip $\{z:\ |{\rm Im}z|<\rho\}$.

Follow the notations in section \ref{sec:Transfer}. We list the corresponding quasiperiodic versions here again for reader's convenience. The analytic quasiperiodic Jacobi operator on ${\ell}^2(\Z)$ is given by:
\begin{equation}
(H_{v,c}u)_n =c(\theta+n\alpha) u_{n+1} + \bar{c}(\theta+(n-1)\alpha)u_{n-1} +v(\theta+n\alpha) u_n,\ n\in\Z.
\end{equation}
The transfer matrix is given by:
\begin{align*}
A(\theta,E,\alpha) = \frac{1}{c(\theta)} \left (\begin{matrix}  E-v(\theta)   &  -\bar{c}(\theta-\alpha) \\ c(\theta)  &   0 \end{matrix} \right )
\end{align*}
and 
\begin{align*}
A(n;\theta,E,\alpha) =\prod_{j=1}^{n} A\big(\theta+(j-1)\alpha,E,\alpha\big),\ n>0. 
\end{align*}

The spectral singularity in Theorem \ref{thm:speSingular} is reduced to the following lemma about the norm of the transfer matrices, which was proved in \cite{JZ}: 
\begin{lem}[\cite{JZ}, Lemma 3.1]\label{lem:JZ3.1}
	Fix $\alpha\in\R\backslash\Q$ with $\beta=\beta(\alpha)<\infty$ and $\theta\in\T$.  Suppose there is a constant $c>0$ such that for any $E$, there is $\ell_0=\ell(E,\beta,\rho,\theta)$ such that for any $\ell>\ell_0$, the following two estimates hold:
	\begin{equation}\label{eq:sumAk-1}
	\sum_{k=1}^{\ell} \|A(k;\theta,E,\alpha)\|^2\ge \ell^{1+\frac{2c}{\beta}},
	\end{equation}
	and
	\begin{equation}\label{eq:sumAk-2}
	\sum_{k=1}^{\ell} \|A(k;\theta-\alpha,E,-\alpha)\|^2\ge \ell^{1+\frac{2c}{\beta}},
	\end{equation}
	then we have the following upper bound for the spectral dimension defined in (\ref{def:dimspe}) of the spectral measure $\mu=\mu_{\alpha,\theta}$:
	\begin{align}\label{eq:speupper}
	{\rm dim}_{{\rm spe}}(\mu)\le \gamma_0:=\frac{1}{1+c/\beta}<1.
	\end{align}
	
\end{lem}

This is a direct consequence of the subordinate theory (\ref{eq:jl}) and Last-Simon upper bound on the generalized eigenfunction (\ref{eq:ls99}). Actually, in view of Lemma \ref{lem:Mmm}, it is enough to find a $\varphi$ such
that both $m_\varphi$ and $\widetilde{m}_{\pi/2-\varphi}$ are $\gamma$-spectral singular, where $m_\varphi$ and $\widetilde{m}_{\pi/2-\varphi}$ are half line m-function defined in section \ref{sec:mMmu}. The estimate on the half line $m$-function relies on the subordinacy theory Lemma \ref{lem:JLineq}. The quantitative estimates need both an upper bound and a lower bound on the $\ell$-norm of $u^\varphi,v^\varphi$. Lemma \ref{lem:ls99} provides two eigen functions $u^{\varphi}$ and $u^{\varphi,-}$, both obeying the sub-linear growth  as in (\ref{eq:ls99}).   (\ref{eq:sumAk-1}) and (\ref{eq:sumAk-2}) provide the lower bound as required in the subordinacy theory for $m_\varphi$ and $\widetilde{m}_{\pi/2-\varphi}$ respectively, which eventually lead to the spectral singularity. In the rest of this section, we will focus on the proof of (\ref{eq:sumAk-1}) and (\ref{eq:sumAk-2}). We refer readers to \cite{JZ}, section 3 for more details about this lemma and spectral singularity. \\

For a ${\rm GL}(2,\mathbb{C})$ matrix 
$A=\left( \begin{array}{cc}  a & b \\ 
 c & d \\        \end{array} \right)$, we denote by $\|\cdot\|_{HS}$ the Hilbert-Smith norm of $A$:
 \begin{align}\label{def:HSnorm}
 \|A\|_{HS}=\sqrt{|a|^2+|b|^2+|c|^2+|d|^2}.
 \end{align}
In the rest of this section, we write  $\|\cdot\|=\|\cdot\|_{HS}$ for simplicity whenever it is clear.

The key to prove (\ref{eq:sumAk-1}) and (\ref{eq:sumAk-2}) is the following lemma:
\begin{lemma}\label{lem:Deltan}
Assume that $L(E)\ge a>0.$ There are $c_2=c_2(a,S,\rho)>0$, $n_0=n_0(a,\rho)>0$ and a
positive integer $d=d(S,\rho,\|v\|_\rho,\|c\|_\rho)\in\N^+$ such that for $E\in S$  and
 $n>n_0$, there exists an interval $\Delta_n\subset \T$ satisfying the following properties:
\begin{equation}\label{eq:LebDelta}
    {\rm Leb}(\Delta_n)\ge \frac{c_2}{4dn}
\end{equation}
and for any $\theta\in \Delta_n$,
\begin{equation}\label{eq:normDelta}
   \|A(n;\theta,E,\alpha)\|^2_{HS}> e^{n L(E)/8}.
\end{equation}

\end{lemma}
Lemma \ref{lem:Deltan} will be the key ingredient to the proof of spectral singularity, we will return to its proof in the end of this section. 
We will derive (\ref{eq:sumAk-1}) and (\ref{eq:sumAk-2}) from  Lemma \ref{lem:Deltan} and finish the proof of Theorem \ref{thm:speSingular} first.

 Let $q_n$ be given as in the continued fraction approximants to $\alpha$, see (\ref{def:beta}). 
The following lemma about the ergodicity of an irrational rotation can be found e.g.  in \cite{JL2}.
\begin{lemma}[Lemma 9, \cite{JL2}]\label{lem:JLDelta}
Let $\Delta\subset [0,1]$ be an arbitrary segment.  If $|\Delta|>\frac{1}{q_n}$. Then, for any $\theta$; there exists a $j$ in
$\{0,1,\cdots,q_n+q_{n-1}-1\}$ such that $\theta+j\alpha\in\Delta$.
\end{lemma}

Combine Lemma \ref{lem:Deltan} with Lemma \ref{lem:JLDelta}, we immediately have the following localization density result:
\begin{lemma}\label{lem:density}
Fix $E\in S,\theta\in\Theta$  and  $\alpha\in\R\backslash\Q$. There is ${n_1}=n_1(E,\rho,\alpha,\theta)$ such that for any $q_n\ge {n_1}$
and  any $m\in\N$, there is $j_m=j_m(\theta)\in\big[2m\,q_n,\,(2m+2)\,q_n\,\big)$ such that
\begin{equation}\label{eq:jN}
    \|A({j_N};\theta,E,\alpha)\|> e^{c_0\,q_n L(E)}
\end{equation}
where $c_0=c_0(a,\rho)$ explicitly depends on $c_2$ and $d$ given in Lemma \ref{lem:Deltan}. 
\end{lemma}
\proof
We fix $E$, $\alpha$ and write $A(n;\theta)=A(n;\theta,E,\alpha)$ for simplicity.
Let $n_0$ be given as in Lemma \ref{lem:Deltan}. Given $q_n$, let 
\begin{equation}\label{def:kn}
    k_n=[\frac{c_2q_n}{4d}]-1\ge \frac{c_2}{5d}q_n\ge n_0,
\end{equation}
provided $q_n$ large, where $c_2$ and $d$ are given as in Lemma \ref{lem:Deltan}.  
By Lemma \ref{lem:Deltan}, there is an interval $\Delta_{k_n}\subset \T$ such that the following hold:
\begin{equation}\label{eq:LebDeltakn}
    {\rm Leb}(\Delta_{k_n})\ge \frac{c_2}{4d{k_n}}> \frac{1}{q_n}
\end{equation}
and 
\begin{equation}\label{eq:normDeltakn}
   \|A(k_n;\theta)\|^2> e^{k_n L(E)/8}>e^{\frac{c_2}{40d}q_n L(E)}, \ \ \forall \theta\in\Delta_{k_n}.
\end{equation}

Fix $\theta$ and $m\in\N$, apply Lemma \ref{lem:JLDelta} to $\Delta_{k_n}$ and $\theta+2m\,q_n$, we have that there exists a $j$ in
$\{0,1,\cdots,q_n+q_{n-1}-1\}$ such that $(\theta+2m\,q_n\alpha)+j\alpha\in\Delta_{k_n}$. By (\ref{eq:normDeltakn}), we have 
\begin{align}\label{eq:Aknexp}
\|A({k_n};\theta+2m\,q_n\alpha+j\alpha)\|> e^{4c_0\,q_n L(E)},
\end{align}
where $c_0=\frac{c_2}{320d}$.

It is easy to check that 
\begin{align}
A({2m\,q_n+j+k_n};\theta)=A({k_n};\theta+2m\,q_n\alpha+j\alpha)\,A({2m\,q_n+j};\theta).
\end{align}
By (\ref{eq:Aknexp}), we have that either
\begin{align}
\|A^{-1}({2m\,q_n+j};\theta)\|\ge e^{2c_0\,q_n L(E)}    \label{eq:Aknj}\\
\ \ \textrm{or} \ \  \|A(2m\,q_n+j+k_n;\theta)\|\ge e^{2c_0\,q_n L(E)}.  \label{eq:Aknj+kn}
\end{align}
Direct computation shows that
\begin{align}
\|A^{-1}({2m\,q_n+j};\theta)\|&=\frac{1}{|{\rm det }A({2m\,q_n+j};\theta)|}\|A({2m\,q_n+j};\theta)\| \\
&= \frac{|c(\theta+({2m\,q_n+j})\alpha)|}{|c(\theta)|}\|A({2m\,q_n+j};\theta)\| \\
&\le \frac{\|c\|_{\infty}}{|c(\theta)|}\|A({2m\,q_n+j};\theta)\|
\end{align}
Suppose (\ref{eq:Aknj}) holds,
then 
\begin{align}
\|A({2m\,q_n+j};\theta)\|\ge \frac{|c(\theta)|}{\|c\|_{\infty}}e^{2c_0\,q_n L(E)} \ge e^{c_0\,q_n L(E)}
\end{align}
provided
\begin{align}
e^{c_0\,q_n L(E)}\ge \frac{\|c\|_{\infty}}{|c(\theta)|}.
\end{align}
Let $j_m$ be $2m\,q_n+j$ or $ 2m\,q_n+j+k_n$, for which $j_N$ satisfies (\ref{eq:jN}).
Clearly, by the choice of $j,k_n$, $j_m(\theta)\in\big[2m\,q_n,\,(2m+2)\,q_n\,\big)$ for all $m\in\N$ and 
\begin{align}\label{eq:sing-n1}
q_n\ge n_1:=\max\big\{   \frac{5dn_0}{c_2},\ \                \frac{\ln \frac{\|c\|_{\infty}}{|c(\theta)|} }{c_0L(E)}               \big \}.
\end{align}
Note that if $m=j=0$ in (\ref{eq:Aknexp}), we pick $j_0=k_n\ge 1$. So $j_0\in \big[1,\,2q_n\,\big)$.
\qed \\

With the above localization density lemma, we can complete the proof of Theorem \ref{thm:speSingular} by checking (\ref{eq:sumAk-1}) and  (\ref{eq:sumAk-2}) in Lemma \ref{lem:JZ3.1} for a.e. $\theta\in\T$.\\

{\noindent \bf Proof of Theorem \ref{thm:speSingular}.}

For any $\ell\in\N$, there is $q_n$ such that, $l\in[2q_n,2q_{n+1}).$ Let $\ell=2Nq_n+r,$ 
where $0\le r<2q_n$, $1\le N<\frac{q_{n+1}}{q_n}$. Let $n_1$ be given as in (\ref{eq:sing-n1}). It is easy to check that 
$q_n\ge n_1$ provided 
\begin{align}
\ell\ge 2e^{2n_1\beta(\alpha)}.
\end{align}

 Now apply Lemma \ref{lem:density} to $q_n$ and $0\le m\le N-1$. There are $j_m\in\big[2m\,q_n,\, (2m+2)\,q_n\big )\subset[0,\ell]$ such that 
$     \|A({j_m};\theta,E,\alpha)\|> e^{c_0\,q_n L(E)}$.  Therefore,
\begin{eqnarray}
  \sum_{k=1}^{\ell} \|A(k;\theta,E,\alpha)\|^2 &\ge & \sum_{m=0}^{N-1}  \|A({j_m};\theta,E,\alpha)\|^2 \nonumber\\
   &\ge&N\,e^{2c_0\,q_n L(E)}. \label{eq:sumAk1-l}
\end{eqnarray}
Clearly,  $\ell=2Nq_n+r<4Nq_n$. By (\ref{eq:sumAk1-l}), we have
$$\sum_{k=1}^{\ell} \|A_k(\theta)\|^2\ge \frac{\ell}{4q_n}e^{2c_0\,q_n L(E)}\ge \ell\;e^{c_0\,q_n L(E)} \ge \ell\;e^{c_0a\, q_n}$$
provided $e^{c_0\,q_n L(E)}\ge 4q_n$. Then for sufficiently large $\ell$ such that $\frac{\ln q_{n+1}}{q_n} <2\beta$, we have

\begin{align}
\sum_{k=1}^{\ell} \|A_k(\theta)\|^2\ge \ell\, q_{n+1}^{\frac{c_0a}{2\beta}}
\ge \ell\cdot \big(\frac{\ell}{2}\big)^{\frac{c_0a}{2\beta}}\ge  \ell\cdot \ell^{\frac{c_0a}{4\beta}} =:\ell^{1+\frac{2c}{\beta}}, \label{eq:pfsumAk}
\end{align}
provided $\ell\ge 4$,  where $c=\frac{1}{8}c_0a$. This proved (\ref{eq:sumAk-1}).

For the same $\theta$ and $E$, repeat the above procedure for
$A(n;\theta-\alpha,E,-\alpha)$. We have a sequence of positive integers $\widetilde{j}_m=\widetilde{j}_m(\theta-\alpha)\in\big[2mq_n,\,2(m+1)q_n\,)$ for any $N\in\N$ and $q_n\ge n_1(E,\rho,-\alpha,\theta-\alpha)$ such that
\begin{equation}
    \|A(\widetilde j_m;\theta-\alpha,E-\alpha,E)\|> e^{c_0\,q_n L(E)}. 
\end{equation}
Note that $c_0=c_0(a,\rho)$ does not depend on $\theta-\alpha$ and is the same as in (\ref{eq:jN}) and  (\ref{eq:pfsumAk}). The same reasoning proves (\ref{eq:sumAk-2}). 

Then by Lemma \ref{lem:JZ3.1}, we have for all $\theta\in\Theta$ and $\beta(\alpha)<\infty$, ${\rm dim}_{\rm spe}(\mu_{\alpha,\theta})<\frac{1}{1+c/\beta}<1$, which completes the proof of Theorem \ref{thm:speSingular}.
\qed \\

In the rest of the section, we focus on the proof of Lemma \ref{lem:Deltan}. In \cite{JZ}, the authors proved the analytic $SL(2, \R)$ version of this lemma. 
One advantage for Shr\"odinger case is the H-S norm $\|A(n;\theta)\|_{HS}^2$  is a real analytic function which can be approximated by trigonometric functions in some uniform sense. For ${\rm GL}(2,\mathbb{C})$ case, the HS norm of the transfer matrices are  meromorphic functions. We need finer decomposition to deal with the poles. 

Fix $E,\alpha$, for $n\in\N^+,\theta\in\T$, let 
\begin{equation}
F_n(\theta)=\|A(n;\theta,E,\alpha)\|_{HS}^2 \label{def:fn}
\end{equation}
be defined as in (\ref{def:HSnorm}). 
We have the following decomposition of $F_n(\theta)$:
\begin{lem}\label{lem:FnDecomp}
For any $E$ and $n\in\N^+$, there are positive functions $f_n(\theta)$ and $g_n(\theta)$ such that 
\begin{align}
&F_n(\theta)=\frac{f_n(\theta)}{g_n(\theta)},  \label{eq:fnDecomp-1} \\ 
& \inf_{n} \frac{1}{n}\int \ln g_n(\theta)\, {\rm d} \theta=0, \qquad  \inf_{n} \frac{1}{n} \int \ln f_n(\theta)\, {\rm d} \theta=2L(E)\label{eq:fnDecomp-2}.
\end{align}
For any $ \varepsilon>0$, there exists $n_0=n_0(\varepsilon)\in\N $ such that for any $n>n_0$  and any $\theta\in\T$, 
\begin{align}
 0< g_n(\theta) <e^{\varepsilon n}  \label{eq:fnDecomp-3}.
\end{align}
Furthermore, for $E$ in a compact set $S$, there are $n_1=n_1(\rho)>0$ and $d=d(S,\rho,\|v\|_\rho,\|c\|_\rho)>0$ such that for any $n>n_1$, 
there are two functions $P_n(\theta),R_n(\theta)$ satisfying the following decomposition:
\begin{align}
& f_n(\theta)=P_n(\theta)+R_n(\theta),   \label{eq:fnDecomp-4}\\
& |R_n(\theta)|<1,   \label{eq:fnDecomp-5}\\
& P_n(\theta)=\sum_{|k|\le d\cdot n}\widehat f_n(k)e^{2\pi {\rm i}k\theta} \label{eq:fnDecomp-6},
\end{align}
where $\widehat f_n(k)$ is the $k$-th Fourier coefficient of $f_n(\theta)$.
\end{lem}

\proof
Follow the notations in (\ref{def:AnDn}), let 
\begin{align*}\label{eq:temp-2}
A(\theta,E)=\frac{1}{c(\theta)}D(\theta,E),\ D(\theta,E)=\left (\begin{matrix}  E-v(\theta)   &  -\bar{c}({\theta}-\alpha) \\ {c}({\theta})  &   0 \end{matrix} \right )
\end{align*}
and
\begin{align}
&A(n;\theta,E)=\frac{1}{c(n;\theta)}D(n;\theta,E),\  \  {\rm where}\\
&c(n;\theta)=\prod ^n_{j=1}c\big(\theta+(j-1)\alpha\,\big), \ \  
D(n;\theta,E)=\prod ^n_{j=1}D\big(\theta+(j-1)\alpha,E\big)=
\left (\begin{matrix}  D^{1}(\theta)   &  D^{2}(\theta) \\ D^{3}(\theta)  &   D^{4}(\theta) \end{matrix} \right ) \label{eq:D1234}
\end{align} 

Without loss of generality, we assume $\int_\T\ln |c(\theta)|{\rm d}\theta=0$. Otherwise, the argument simply differs by a constant factor. See Remark \ref{rmk:lnc=0} after the proof.

Let $g_n(\theta):=\big|c(n;\theta)\big|^2$ and $f_n(\theta):=\|D(n;\theta,E)\|^2_{HS}$. Clearly, 
\begin{align}
&F_n(\theta)=\|A(n;\theta,E,\alpha)\|_{HS}^2=\frac{f_n(\theta)}{g_n(\theta)} \label{eq:fndecomp-pf}\\
&f_n(\theta)=\|D(n;\theta,E)\|^2_{HS}=|D^{1}(\theta)|^2+|D^{2}(\theta)|^2+|D^{3}(\theta)|^2+|D^{4}(\theta)|^2. \label{eq:DnHSnorm}
\end{align} 
Birkhoff Ergodic Theory implies that for any irrational $\alpha$,
\begin{align}
\lim_n\int_\T\frac{1}{n}\ln g_n(\theta)\, {\rm d} \theta
&=\inf_{n} \frac{1}{n}\int \ln g_n(\theta)\, {\rm d} \theta   \nonumber \\
&=\lim_n\int_\T \frac{1}{n}\sum_{j=1}^n\ln |c\big(\theta+(j-1)\alpha\,\big)|^2\, {\rm d} \theta
=\int_\T\ln|c(\theta)|^2\,{\rm d}\theta=0
\end{align} 
In view of (\ref{eq:fndecomp-pf}) and the definition of Lyapunov exponent (\ref{def:Lyp}), we have 
\begin{align}
\inf_{n}  \int \frac{1}{n}\ln f_n(\theta) \, {\rm d} \theta
&= \inf_{n}  \int \frac{1}{n}\ln\Big(g_n(\theta)\,\|A(n;\theta)\|_{HS}^2\Big)\, {\rm d} \theta \nonumber\\
&=\inf_{n}  \int \frac{1}{n}\ln g_n(\theta)\, {\rm d} \theta+\inf_{n} \int  \frac{1}{n}\ln \|A(n;\theta)\|_{HS}^2\, {\rm d} \theta \nonumber  \\
&=2L(E). 
\end{align} 
Note  $c(\theta)$ is continuous in $\theta$, by (\ref{eq:limsup-a}),for any $\varepsilon>0$, there is $n_1=n_1(\varepsilon)$ such that for any $n>n_1$ and any $\theta\in\T$, we have the following upper semicontinuity (uniform in $\theta$):
\begin{align}
 \frac{1}{n}\ln g_n(\theta)\le \int_\T\ln|c(\theta)|^2\,{\rm d}\theta+\varepsilon=\varepsilon. \label{eq:uniform-g}
\end{align}
This gives $g_n(\theta)\le e^{\varepsilon n}$ and finishes the proof of (\ref{eq:fnDecomp-1})-(\ref{eq:fnDecomp-3}).\\

The further decomposition of $f_n(\theta)$ into $P_n$ and $R_n$ follows the strategy in \cite{JZ}. Note that $v(\theta)$ and $c(\theta)$ are both analytic with bounded extension to the strip 
$\{z:|Imz|<\rho\}$. In view of (\ref{eq:D1234}), all $D^{i}(\theta),i=1,2,3,4$ have analytic
extension to the strip $\{z:|Imz|<\rho\}$.  For compact $S$, there is $C_1=C_1(S,\rho,\|v\|_\rho,\|c\|_\rho)$ such that
\begin{equation}\label{eq:Dnrho}
   \|D^{i}\|_{\rho}:=\sup_{|Imz|<\rho} \Big|D^{i}(z)\Big|<\sup_{|Imz|<\rho}\|D_n(z)\|^2_{HS}<e^{C_1n}, \ \ E\in S, \ \ i=1,2,3,4.
\end{equation}
 Consider the Fourier expansion of the periodic-1 functions $D^{i}(\theta)$:
\begin{equation}\label{eq:DnSum}
    D^{i}(\theta)=\sum_{k\in\Z}\widehat D^i(k)e^{2\pi ik\theta}, \ \ i=1,2,3,4 
\end{equation}
The Fourier coefficients of $ D^{i}(\theta)$ has exponential decay as 
\begin{align}
  |\widehat D^i(k)|<\|D^i\|_{\rho}\cdot e^{-2\pi\rho|k|}<e^{C_1n}\cdot e^{-2\pi\rho|k|}, \ \forall k\in\Z, \ \ \ \ i=1,2,3,4.  \label{eq:DikDecay}
\end{align}
Combine \begin{equation}\label{eq:DbarDSum}
    |D^{i}(\theta)|^2=\big(\sum_{k\in\Z}\widehat D^i(k)e^{2\pi ik\theta}\big)\,\big(\sum_{k\in\Z} \overline {\widehat {D^i}(k)}e^{-2\pi ik\theta}\big)
\end{equation}
with (\ref{eq:DikDecay}), it is easy to check that the Fourier coefficients of $ |D^{i}(\theta)|^2$ has exponential decay as:
\begin{align}
 \big |\,\widehat {|D^i|^2(\cdot)}(k)\,\big|<e^{C_1n}\cdot e^{-\pi\rho|k|}, \ \forall k\in\Z, \ \ \ \ i=1,2,3,4 . \label{eq:DbarDDecay}
\end{align}

Let $ f_n(\theta)$ be given as in (\ref{eq:DnHSnorm}). Consider the Fourier expansion of $ f_n(\theta)$:
\begin{equation}\label{eq:fnSum}
    f_n(\theta)=\sum_{k\in\Z}\widehat f_n(k)e^{2\pi ik\theta}.
\end{equation}
By (\ref{eq:DnHSnorm}) and (\ref{eq:DbarDDecay}), clearly, $\widehat f_n(k)$ has the same exponential decay in $|k|$:
\begin{align}
  |\widehat f_n(k)|<4e^{C_1n}\cdot e^{-\pi\rho|k|}, \ \forall k\in\Z.
\end{align}

Pick
\begin{equation}\label{def:d}
d=\left[\frac{C_1}{\pi\rho}\right]+2.
\end{equation}
We split $f_n(\theta)$ into two parts:
$$
  f_n(\theta)=P_n(\theta)+R_n(\theta),\
  P_n(\theta)=\sum_{|k|\le d\cdot n}\widehat f_n(k)e^{2\pi ik\theta},\ R_n(\theta)= \sum_{|k|>d\cdot n}\widehat f_n(k)e^{2\pi ik\theta}.
$$
For any $\theta\in\T$,
\begin{eqnarray*}
  |R_n(\theta)| \le  \sum_{|k|>d\cdot n}|{\widehat f_n(k)}| &\le &  \sum_{|k|>d\cdot n}4e^{C_1n}\cdot e^{-\pi\rho|k|} \\
   &\le& \frac{8}{1-e^{-\pi\rho}}e^{C_1n}e^{-\pi\rho dn}\\
   &\le & \frac{8}{1-e^{-\pi\rho}}e^{-(\pi\rho\, d-C_1)n}.
\end{eqnarray*}
By the choice of  $d$ in (\ref{def:d}), we have $\pi\rho d>C_1+\pi\rho$. Then for any $\theta\in\T$,
\begin{equation}\label{eq:rn}
    |R_n(\theta)|\le \frac{8}{1-e^{-\pi\rho}}e^{-\pi\rho\,n}<1,
\end{equation}
provided $n>n_2(\rho):=(\pi\rho)^{-1}\ln(\frac{8}{1-e^{-\pi\rho}})$. 
This finishes the proof of (\ref{eq:fnDecomp-4})-(\ref{eq:fnDecomp-6}).
\qed\\

\begin{rem}\label{rmk:lnc=0}
	Suppose $b=\int_\T\ln|c(\theta)|{\rm d}\theta\neq0$. In (\ref{eq:temp-2}), we set
	\begin{align}
	A(\theta,E)=\frac{1}{\widetilde c(\theta)}\,\widetilde D(\theta,E),\ {\rm where}\ \ 
	\widetilde c(\theta)=e^{-b}c(\theta),\ \widetilde D(\theta,E)=e^{-b}D(\theta,E).
	\end{align}
	Clearly, 
	\begin{align}
\int_\T\ln	|\widetilde c(\theta)|{\rm d}\theta=0,\ 
\lim_n \int_\T\frac{1}{n}\ln  \|\widetilde D(n;\theta,E)\|{\rm d }\theta=L(E).
	\end{align}
Let $g_n(\theta):=\big|\widetilde c(n;\theta)\big|^2$ and $f_n(\theta):=\|\widetilde D(n;\theta,E)\|^2_{HS}$. The rest of the decomposition are exactly the same. 
\end{rem}


Combine Lemma \ref{lem:FnDecomp} with the positive assumption on Lyapunov exponent, we can now finish {\noindent \bf The proof of Lemma \ref{lem:Deltan}:}
Assume that the Lyapunov exponent $L(E)\ge a>0$ for $E\in S$. Pick $\varepsilon=a/8$. Let $n_1=n_1(\varepsilon)$ and $n_2=n_2(\rho)$ be given as in Lemma \ref{lem:FnDecomp}. Then for all $n>\max\{n_1,n_2\}$, we have  $g_n(\theta),f_n(\theta),P_n(\theta)$ and $R_n(\theta)$ as in Lemma \ref{lem:FnDecomp}, satisfying (\ref{eq:fnDecomp-1})-(\ref{eq:fnDecomp-6}).  Denote
\begin{eqnarray*}
  \Theta_n^1 &=&\{\theta: F_n(\theta)> e^{n L(E)/8} \}, \\
  \Theta_n^2 &=&\{\theta: P_n(\theta)> e^{n L(E)/3} \}, \\
  \Theta_n^3 &=&\{\theta: f_n(\theta)> e^{n L(E)/2} \}.
\end{eqnarray*}
Let $n_3:=4a^{-1}$. Then for all $n>n_3$, we have  $
	e^{nL(E)}>e^{na}>e^4>50$. 
By using the fact $x^{1/2}-x^{1/3}>x^{1/3}-x^{1/4}>1$ for all $x>50$, it is easy to check that for $n>n_3$, 
\begin{align}\label{eq:nL>1}
e^{n L(E)/2}-e^{n L(E)/3}>e^{n L(E)/3}-e^{n L(E)/4}>1.
\end{align}
Assume that $f_n(\theta)> e^{n L(E)/2}$. 
By (\ref{eq:fnDecomp-4}) and \eqref{eq:nL>1}, we have for $n>n_3$, 
$$P_n(\theta)>f_n(\theta)-|R_n(\theta)|>e^{n L(E)/2}-1>e^{n L(E)/3}. $$
Then
$$f_n(\theta)>P_n(\theta)-|R_n(\theta)|>e^{n L(E)/3}-1>
e^{n L(E)/4}.$$
In view of (\ref{eq:fnDecomp-1}) and (\ref{eq:fnDecomp-3}), we have then for $n>\max\{n_1,n_3\}$,
$$F_n(\theta)=\frac{f_n(\theta)}{g_n(\theta)}>
\frac{e^{n L(E)/4}}{e^{n \varepsilon}}>\frac{e^{n L(E)/4}}{e^{n L(E)/8}}=
e^{n L(E)/8}.$$
Therefore, we have for $n>n_0:=\max\{n_1,n_2,n_3\}$, 
\begin{equation}\label{eq:theta}
   \Theta_n^3\subseteq\Theta_n^2\subseteq\Theta_n^1.
\end{equation}
Meanwhile, by (\ref{eq:fnDecomp-2}), 
\begin{eqnarray*}
  2nL(E) &\le& \int_{\T}\ln f_n(\theta){\rm d}\theta \\
   &\le&  {\rm Leb}(\Theta_n^3) \ln \|f_n\|_{\rho}+ \big(1-{\rm Leb}(\Theta_n^3)\big) \ln e^{n L(E)/2} \\
   &\le&  {\rm Leb}(\Theta_n^3) \cdot C_1n+ \big(1-{\rm Leb}(\Theta_n^3)\big) \cdot n L(E)/2.
\end{eqnarray*}
This implies ${\rm Leb}(\Theta_n^3)\ge \frac{3L(E)}{2C_1-L(E)}$. Note that 
$L(E)\ge a>0,E\in S$, we have
\begin{equation}\label{eq:c2}
    {\rm Leb}(\Theta_n^3)\ge \frac{3a}{2C_1-a}=:c_2(a,S,\rho)>0. 
\end{equation}
In view of (\ref{eq:theta}), we have for $n>n_0$, 
\begin{equation}\label{ce2}
    {\rm Leb}(\Theta_n^2)
\ge c_2(a,S,\rho)>0. 
\end{equation}
By (\ref{eq:fnDecomp-6}), $P_n(\theta)$ is a trigonometric polynomial of degree (at most)
$2dn$, where $d$ is given by (\ref{def:d}) in Lemma \ref{lem:FnDecomp}. The set $\Theta_n^2$ consists of no more than $4dn$ intervals. Therefore, there exists a
segment, $\Delta_n\subset \Theta_n^2\subset \Theta_n^1$, with ${\rm Leb}(\Delta_n)>\frac{c_2}{4dn}$. 
For any $n>n_0$ and $\theta\in \Delta_n\subset \Theta_n^1$, 
$$\|A_n(\theta)\|^2_{HS}=F_n(\theta)> e^{n L(E)/8}$$ and
$${\rm Leb}(\Delta_n)>\frac{c_2}{4dn},$$
as claimed. \qed

\section{The Extended Harper's model: proof of Corollary \ref{cor:EHM} } \label{sec:EHM}
Recall the extended Harper's model (EHM) defined in (\ref{def:EHM}) as:
\begin{equation}\label{eq:EHM}
(H_{\lambda,\alpha,\theta}u)_n=c_{\lambda}(\theta+n\alpha)u_{n+1}+\bar{c}_{\lambda}\big(\theta+(n-1)\alpha\big)u_{n-1}+2\cos{2\pi (\theta+n\alpha)} u_n,
\end{equation}
where
\begin{align} \label{eq:clambda}
c_{\lambda}(\theta)=\lambda_1 e^{-2\pi i(\theta+\frac{\alpha}{2})}+\lambda_2 +\lambda_3 e^{2\pi i(\theta+\frac{\alpha}{2})}.
\end{align}

By some earlier work \cite{JMarx12CMP}, we consider the following partitioning of the parameter space into the following three regions:
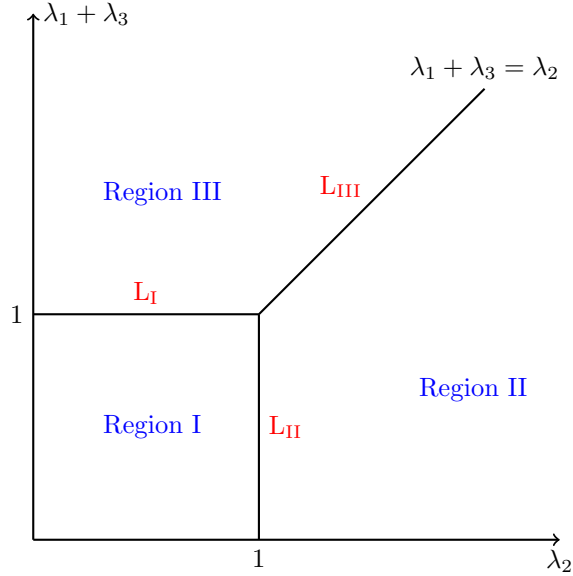
\begin{center}
	\begin{tikzpicture}[thick, scale=1]
	\draw[->] (-10,-1) -- (-3,-1) node[below] {$\lambda_2$};
	\draw[->] (-10,-1) -- (-10,6) node[right] {$\lambda_1+\lambda_3$};
	\draw [ ] plot [smooth] coordinates { (- 7, 2) (-4, 5) };
	\draw [ ] plot [smooth] coordinates { (- 7, 2) (-7,-1) };
	\draw [ ] plot [smooth] coordinates { (-10, 2) (-7, 2) };
	
	\draw(-  4,   5) node [above] {$\lambda_1+\lambda_3=\lambda_2$};
	\draw(-  7,  -1) node [below] {$1$};
	\draw(- 10,   2) node [left]  {$1$};
	\draw(-9.2, 0.5) node [color=blue][right] {Region I};
	\draw(-  5,   1) node [color=blue][right] {Region II};
	\draw(-9.2, 3.6) node [color=blue][right] {Region III};
	
	\draw(-  7, 0.5) node [color=red][right] {$\mathrm{L}_{\mathrm{II}}$};
	\draw(-8.5,   2) node [color=red][above] {$\mathrm{L}_{\mathrm{I}}$};
	\draw(-5.5, 3.7) node [color=red][left] {$\mathrm{L}_{\mathrm{III}}$};
	
	\end{tikzpicture}
\end{center}

\begin{description}
	\item[Region I] $0 \leq \lambda_{1}+\lambda_{3} \leq 1, ~0 < \lambda_{2} \leq 1$ ~\mbox{,}
	\item[Region II] $\max\{\lambda_{1}+\lambda_{3},1\} \leq \lambda_{2}$, $\lambda_1+\lambda_3>0$~\mbox{,}
	\item[Region III] $\max\{1,\lambda_{2}\} \leq \lambda_{1}+\lambda_{3}$, $\lambda_2>0$ ~\mbox{.}
\end{description}

Let $L(E,\lambda)$ be the {\it{Lyapunov exponent}} of the extended Harper's model, defined as in (\ref{def:Lyp}).  The main achievement of \cite{JMarx12CMP} is to prove the following explicit formula of $L(E,\lambda)$, valid for all $\lambda$ and all irrational $\alpha$:
\begin{thm}[\cite{JMarx12CMP}] \label{thm:JM12CMP}
	Fix an irrational frequency $\alpha$. Then $L(E,\lambda)$ restricted to the spectrum is zero within both region II and III. In region I it is given by the formula on the spectrum,
	\begin{equation} \label{def:EHMLyp}
L(E,\lambda)=	\begin{cases}\ln \left(  \dfrac{1+\sqrt{1 - 4\lambda_{1} \lambda_{3}}}{2 \lambda_{1}}\right) & \mbox{, if} ~\lambda_{1} \geq \lambda_{3}, ~\lambda_{2} \leq \lambda_{3} + \lambda_{1} ~\mbox{,}  \\
	\ln \left(  \dfrac{1+\sqrt{1 - 4\lambda_{1} \lambda_{3}}}{2 \lambda_{3}}\right) & \mbox{, if} ~\lambda_{3} \geq \lambda_{1}, ~\lambda_{2} \leq \lambda_{3} + \lambda_{1} ~\mbox{,}  \\
	\ln \left(  \dfrac{1+\sqrt{1 - 4\lambda_{1} \lambda_{3}}}{\lambda_{2} + \sqrt{\lambda_{2}^{2} - 4 \lambda_{1} \lambda_{3}}}\right) & ~\mbox{, if} ~\lambda_{2} \geq \lambda_{3} + \lambda_{1} ~\mbox{.} 
	\end{cases}
	\end{equation}
\end{thm}

Denote by Region ${\rm I}^\circ$,Region ${\rm II}^\circ$, Region ${\rm III}^\circ$ the interior of Region I,II,III respectively. A complete understanding of the spectral properties of the extended Harper's model for a.e. $\theta$ has been established in \cite{JKS,HJ,AJM,H17}. 
We collect the spectral decomposition results in these papers as the follow theorem for reader's convenience. Follow the notations in Corollary \ref{cor:EHM}, denote the three parameter regions of $\lambda=(\lambda_1,\lambda_2,\lambda_3)\in\R^3$ by:
\begin{align*}
&{\cal R}_1=\left\{\lambda\in\R^3:\, 0<\lambda_1+\lambda_3<1,0<\lambda_2<1 \right\}.\\
&{\cal R}_2=\left\{\lambda\in\R^3:\, \lambda_2>\max\{\lambda_1+\lambda_3,1\},\lambda_1+\lambda_3\ge0 \ \ 
{\rm or}\ \ \lambda_1+\lambda_3>\max\{\lambda_2,1\},\lambda_1\neq\lambda_3,\lambda_2>0 \right\} .\\
&{\cal R}_3=\{\lambda\in\R^3:\, 0\le \lambda_1+\lambda_3\le 1,\lambda_2=1 \ \  {\rm or}\ \ \lambda_1+\lambda_3\ge \max\{\lambda_2,1\},\lambda_1=\lambda_3,\lambda_2>0 \} .
\end{align*}
\begin{thm}[\cite{JKS,HJ,AJM,H17}]\label{thm:EHMspe}The following Lebesgue decomposition of the spectrum of $H_{\lambda,\alpha,\theta}$ holds for a.e. $\theta$.
	\begin{itemize}
		\item For  $\lambda\in{\cal R_1}$, if $\beta(\alpha)<L(E,\lambda)$, then $H_{\lambda,\alpha,\theta}$ has pure point spectrum. If $\beta(\alpha)>L(E,\lambda)$, then $H_{\lambda,\alpha,\theta}$ has purely singular continuous spectrum.
		\item For $\lambda\in{\cal R_2}$ and all irrational $\alpha$, $H_{\lambda,\alpha,\theta}$ has purely absolutely continuous spectrum.
		\item For $\lambda\in{\cal R_3}$ and all irrational $\alpha$, $H_{\lambda,\alpha,\theta}$ has purely singular continuous spectrum.
	\end{itemize}
\end{thm}
Now we are in the place to analyze the spectral dimension of EHM in each region. 

Clearly, Region ${\rm I}^\circ={\cal R}_1$. In view of (\ref{def:EHMLyp}), it is easy to check that $L(E,\lambda)>0$ on ${\cal R}_1$ for all $\alpha$ and $E$. Therefore, by Theorem \ref{thm:dimspe=1}, we have part (1) of Corollary \ref{cor:EHM}.

Next, consider Region ${\cal R}_2$. Theorem \ref{thm:EHMspe} shows that $H_{\lambda,\alpha,\theta}$ has purely a.c. spectrum in region  ${\cal R}_2$ for all $\alpha$ and a.e. $\theta$. In view of Definition \ref{def:dimspe}, absolutely continuous measure has full spectral dimension\footnote{Actually, it is well known that a.c. measure has full dimension for most commonly used fractal dimensions, e.g. Hausdorff/packing dimension etc. See more background knowledge about fractal dimension in e.g. \cite{Fal}}. This gives part (2) of Corollary \ref{cor:EHM}.

Part (3) (region ${\cal R}_3$) is the only place requires extra work. 
 By Theorem \ref{thm:JM12CMP} and Theorem \ref{thm:EHMspe}, in region  ${\cal R}_3$,  $L(E,\lambda)=0$ on the spectrum and $H_{\lambda,\alpha,\theta}$ does not have a.c. spectrum. Lack of positivity of Lyapunov exponent, we do not have the spectral singularity and the upper bound provided by Theprem \ref{thm:speSingular}. While the lower bound from Theorem \ref{thm:speconti} still holds. Moreover, in view of Lemma \ref{lem:Furman} and Corollary \ref{cor:c-analytic-bound}, we can obtain arbitrarily small exponential growth of the transfer matrix. This allows us to obatian the increased range of $\beta(\alpha)$ in the critical region in part (3).

 Recall the notations of the tranfer matrix in (\ref{def:Dnm}) and (\ref{def:wnm}) for EHM:
let
\begin{align*}
A^{\lambda}(\theta,E,\alpha) = \frac{1}{c_{\lambda}(\theta)}D^{\lambda}(\theta,E,\alpha),\ \ 
D^{\lambda}(\theta,E,\alpha)=\left (\begin{matrix}  E-v(\theta)   &  -\bar{c_{\lambda}}(\theta-\alpha) \\ c_{\lambda}(\theta)  &   0 \end{matrix} \right ).
\end{align*}
For $n>0,m\in\Z$,
\begin{align}
&A^{\lambda}(n,m;\theta) =\prod_{j=m}^{m+n-1} A^{\lambda}\big(\theta+j\alpha\big), \\
&D^{\lambda}(n,m;\theta) =\prod_{j=m}^{m+n-1} D^{\lambda}\big(\theta+j\alpha\big),\qquad
c_{\lambda}(n,m;\theta) =\prod_{j=m}^{m+n-1} c_{\lambda}\big(\theta+j\alpha\big). \label{eq:DcEHM}
\end{align}

It is easy to check that 
\begin{align}
L(E,\lambda)=L(D^{\lambda})-\int_\T\ln|c_\lambda(\theta)|{\rm d}\theta.
\end{align}
Note that
\begin{align}
b_\lambda:=\int_\T\ln|c_\lambda(\theta)|{\rm d}\theta
\end{align}
 is not necessarily zero in region ${\cal R}_3$. Suppose not, consider the rescaling trick in Remark \ref{rmk:lnc=0}. Set
 \begin{align}
 \widetilde c_{\lambda}(\theta)=e^{-b_\lambda}c_{\lambda}(\theta),\ \widetilde D^\lambda(\theta,E)=e^{-b_\lambda}D^\lambda(\theta,E).
 \end{align}
 Clearly, in Region X,
 \begin{align}
 \int_\T\ln	|\widetilde c(\theta)|=0,\ 
 L(\widetilde D^\lambda)=L(D^{\lambda})-b_\lambda=L(E,\lambda)=0.
 \end{align}
 Let $\widetilde D^{\lambda}(n,m;\theta),\widetilde c_{\lambda}(n,m;\theta)$ be defined the same way as in (\ref{eq:DcEHM}). 
For irrational $\alpha$, let $\beta(\alpha)$ and $q_n$\footnote{We still denote the subsequence reaching the $\limsup$ by $q_n$. } be defined as in (\ref{def:beta}).  Now assume $\beta(\alpha)>0$, let $\widetilde\beta=\min\{\beta(\alpha)/3,1\}$. 
It was proved in \cite{JMarx12CMP} that $L(E,\alpha)$ is continuous in $E$ for irrational $\alpha$. In view of Lemma \ref{lem:Furman}, the $\limsup$ is uniform in both $\theta$ and $E$.  Therefore, for any $\delta>0$, there is $n_0=n_0(\delta, \widetilde\beta)$ such that for any $n>n_0$, $m\in\Z$, $\theta\in\T$ and $E\in \sigma(H_{\lambda,\alpha,\theta})$,  
\begin{align}
& \|\widetilde D^{\lambda}(n,m;\theta)\| \leq e^{\delta^2\,\widetilde\beta\, n},\label{eq:EHMD} \\
& |\widetilde c_{\lambda}(n,m;\theta)| \leq  e^{\delta^2\,\widetilde\beta\, n}\label{eq:EHMc} .
\end{align}
Note that in the proof Theorem \ref{thm:trace}, we only need to consider the above upper bound for $E$ restricted in the spectrum.  
By Corollary \ref{cor:c-analytic-bound},  for a.e. $\theta$ and $q_n$ large, 
 \begin{align}
 \min_{|m|\le e^{\delta\beta\,q_n}}|\widetilde c_{\lambda}(q_n,m;\theta)|>e^{-6\delta^2 \widetilde\beta\,q_n}. \label{eq:EHMclower}
 \end{align}
Combing (\ref{eq:EHMD}), \eqref{eq:EHMc} and \eqref{eq:EHMclower}, exact the same computation in section \ref{sec:MoreLambdaBound} shows that for a.e. $\theta$,  $0<\delta<\frac{1}{\sqrt 7}$ and $q_n$ large, 
 \begin{align}
&\min_{|m|\le e^{\delta\beta\,q_n}} |\widetilde c_{\lambda}(r,m;\theta)| \ge e^{-7\delta^2 \widetilde\beta\,q_n}, \ \ 1\le r\le q_n,\\
&\max_{|m|\le e^{\delta\beta\,q_n}}\Big|\frac{\widetilde c_{\lambda}({\theta+(m\pm q_n)\alpha})}{\widetilde c_{\lambda}({\theta+m\alpha})}-1\Big| < e^{-(\beta-7\delta^2 \widetilde\beta) q_n},\\
&
\sup_{E\in \sigma(H_{\lambda,\alpha,\theta})}\|A^{\lambda}(r,m;\theta)\|< e^{8\delta^2 \widetilde\beta q_n},\ \ 0\le r\le q_n,\ |m|\le e^{\delta\beta\,q_n}. 
\end{align} 
 Therefore, we can replace all the $\Lambda$ in the proof Theorem \ref{thm:trace} by $10\delta^2 \widetilde\beta$. Then for any $\beta(\alpha)>0$ and $0<\gamma<1$, (\ref{eq:beta>Lambda}) holds true provided 
\begin{align}
\delta<\frac{1}{6000}(1-\gamma). 
\end{align}
Thereofore, by Lemma \ref{lem:Powerlaw} and Theorem \ref{thm:trace}, for any $\beta(\alpha)>0$, $\gamma<1$ and a.e. $\theta$, $\mu_{\lambda,\alpha,\theta}$ is $\gamma$-spectral continuous. By (\ref{eq:dimspe}), ${\rm dim}_{\rm spe}(\mu_{\lambda,\alpha,\theta})$=1, which completes the proof of part (3) of Corollary \ref{cor:EHM}. \qed
\appendix
\section{Appendix}

\subsection{Proof of (\ref{eq:beta>Lambda}) in Theorem \ref{thm:trace}}\label{A:pflemPower}
We have showed in the first part of Theorem \ref{thm:trace} that if $\beta>260(1+\frac{1}{\delta})\Lambda$, then for  $\mu\ a.e.\ E$, there exists $K(E)\in\N$, for $k\geq K(E)$, we have
\begin{align} \label{eq:app-trace}
|Tr A(q_k;E)| <2-2e^{-60\Lambda q_k}
\end{align}

Now by (\ref{eq:tr-tr}), we have 
\begin{align}\label{eq:A1-1}
|Tr \widetilde A(q_k;E)| <2-2e^{-60\Lambda q_k}+12e^{(-\beta+4\Lambda)q_k}<2-e^{-60\Lambda q_k},
\end{align}
provided $e^{(\beta-64\Lambda)q_k}>12$. Fix $E$ and $q=q_k$ and write $\widetilde A(q_k;E)=\widetilde A(q)$. 
Now apply Lemma (\ref{lem:JZ2.9}) to these $\widetilde A(q)$ satisfying \ref{eq:A1-1}. Note $\widetilde A(q)\in {\rm SL}(2,\R)$, and $|{\rm Tr}\widetilde A(q)|<2$, the eigenvalue $\rho$ of $\widetilde A(q)$ is purely imaginary with modulus 1, i.e.,   $\rho=e^{ i\psi}$,  for some $\psi\in(-\pi,\pi)$. By (\ref{eq:JZ2.9-1}), we have for any $j$,
\begin{equation}\label{eq:Akq6}
\widetilde A^j(q)=\frac{\sin j\psi}{\sin \psi}\cdot\Big(
\widetilde A(q)-\frac{{\rm Tr\,}\widetilde A(q)}{2}\cdot I\Big)+\frac{\cos j\psi}{2}\cdot I, \ \ \psi\in(-\pi,\pi)
\end{equation}
Then $|2\cos\psi|=|{\rm Tr\,}\widetilde A(q)|<2-e^{- 60\Lambda q}$ implies
$|\sin\psi|
>\sqrt{1-(1-\frac{1}{2}e^{- 60\Lambda q})^2}
>e^{- 40\Lambda q}$.
By (\ref{eq:Akq6}) and (\ref{eq:normA}),
\begin{align}\label{eq:A1-normA}
\|\widetilde A^j(q)\|\le 2 e^{ 40\Lambda q}\,\|\widetilde A(q)\|+1\le e^{ 43\Lambda q},
\end{align}
provided $q>q(\Lambda)$. 


Now for any $0<\gamma<1$, let $\xi=\frac{95}{1-\gamma}<e^{\delta \beta q}$ and  
\begin{align}\label{def:N-gamma}
N=[e^{\xi \Lambda q}].
\end{align}
Apply Lemma \ref{lem:JZA1} to $G=\widetilde A(q)$, $G_j=\widetilde A(q,jq+1)$, $j=0,\cdots, N$, by (\ref{eq:tildeAbeta}) and (\ref{eq:A1-normA}),  for all $j\le N$ we have
\begin{align}
\|\widetilde A({jq})-\widetilde A^j(q)\|<e^{\big(-\beta+93\Lambda+2\xi\Lambda\big)q}<e^{-\Lambda q}<1.
\end{align}
provided $\beta>(94+2\xi)\Lambda$. Therefore, by (\ref{eq:A1-normA}), 
$$\|\widetilde A({jq})\|\le \|\widetilde A^j(q)\|+1\le 2 e^{ 40\Lambda q}\,\|\widetilde A(q)\|+2\le e^{43\Lambda q} $$
By (\ref{eq:rexpbeta}), $|r^{-1}(jq,0)|\le 1+Ne^{(-\beta+2\Lambda)q}\le 1+e^{(-\beta+2\Lambda+\xi\Lambda)q}<2$ provided $\beta>3\Lambda+\xi\Lambda$.
Then by (\ref{eq:conjA-A}) and (\ref{eq:normA}), for all $0\le j\le N$ and $1\le r\le q$,
\begin{align}
&\| A({jq})\|\le |r^{-1}(jq,0)|\cdot\|T_{jq}\|\cdot\|\widetilde A({jq})\|\cdot\|T_0^{-1}\|\le 2e^{43\Lambda q}\\
&\| A({jq}+r)\|\le \| A(r,{jq}+1)\|\cdot \| A({jq})\|\le e^{46\Lambda q}
\end{align}
Therefore,
\begin{align}
&\sum_{n=1}^{N q}\|A(n;E)\|^2 \le\sum_{k=0}^{N-1}\sum_{r=1}^{q}\|A({kq+r};E)\|^2
\le N q\, e^{92\Lambda q}
\le e^{(\xi+93)\Lambda q}\\
&\frac{1}{(Nq)^{2-\gamma}}\sum_{n=1}^{N q}\|A(n;E)\|^2\le e^{\big(-(1-\gamma)\xi+94\big)\Lambda q}=e^{-\Lambda q}<1
\end{align}
In conclusion, for any $0<\gamma<1$ and $\mu$ a.e. $E$, we have a sequence $q_k\to \infty$ and $\ell_k=[e^{95(1-\gamma)^{-1}\Lambda q_k}]q_k$ such that 
\begin{align}\label{eq-1}
\sum_{n=1}^{\ell_k}\|A(n;E)\|^2\le \ell_k^{2-\gamma}
\end{align}
provided 
\begin{align}
\beta>(3\xi+\xi/\delta)\Lambda=(285+\frac{95}{\delta})\frac{\Lambda}{1-\gamma}>(94+2\xi+\xi/\delta)\Lambda.
\end{align}
It was proved in \cite{JZ} that (\ref{eq-1}) implies (\ref{eq:uvLk}) directly from the relation (\ref{eq:eigeneq}) and (\ref{eq:subor-l}) and . We omit the proof for this part here. See more details about this direct computation in the proof Lemma 2.1 in \cite{JZ}.  
\qed

\subsection{Proof of Claim \ref{clm:JZclm5}}\label{A:pfclm5}
For any $0<\varepsilon<e^{-29\Lambda q}$, let $\ell=\ell(\varphi,\varepsilon,E),u^{\varphi},v^{\varphi}$ be given as in (\ref{eq:subor-l}). Write $\ell(\varepsilon)=[\ell]+\ell-[\ell]$, and $[\ell]=K(\varepsilon)\cdot q+r(\varepsilon)$, where $0\le r=[\ell]{\rm mod }\ q< q$
and $0\le \ell-[\ell]<1$. Let $X=\left(
\begin{array}{c}
\cos\varphi \\
-\sin\varphi \\
\end{array}
\right)
$ and $\widetilde X=T_0^{-1}X$. Clearly, $\|X\|=\|\widetilde X\|=1$.

We need to show first suppose $K<N_q=[e^{{29\Lambda q}}]$, then for any $\varepsilon<e^{-29\Lambda q}$:
\begin{equation}\label{eq:K-lower}
K>e^{\Lambda q}
\end{equation}

 For any $n\le [\ell]+1$, write $n=kq+r$, where $0\le k\le K, \ 0\le r\le q$.
By (\ref{eq:rexpbeta}), (\ref{eq:normAlinear}) and (\ref{eq:A(kq)-Ak(q)}), we have
$$\| A(kq)\|\le |r^{-1}(kq,0)|\cdot (\|\widetilde A^k(q)\|+1)\le 2\,(6C_1\,k\,e^{2\Lambda q}+1)<k\, e^{3\Lambda q}$$
Then by (\ref{eq:normA}),
$$\|A(kq+r)\, X\|\le \|A(r,kq+1)\|\cdot \| A(kq)\|\cdot \|X\|\le k\, e^{5\Lambda q}$$

Direct computation shows
\begin{eqnarray*}
	\|u^{\varphi}\|_{\ell}^2    \le   \sum_{n=1}^{[\ell]+1}\|A(n)\cdot X\|^2 
	&\le & \sum_{r=1}^{q}\|A({r})\cdot X\|^2
	+\sum_{k=1}^{K}\sum_{r=1}^{q}\|A({kq+r})\cdot X\|^2 \\
		&\le &  q\cdot e^{4\Lambda q}
	+\sum_{k=1}^{K}\sum_{r=1}^{q}k^2e^{10\Lambda q} \\
	&\le& q\cdot e^{4\Lambda q} +
	K^3\,q\, e^{10\Lambda q}\\
	&\le& K^3\,e^{11\Lambda q}
\end{eqnarray*}

Since $\varphi$ is arbitrary, we have  $\|v^{\varphi}\|_{\ell}^2\le K^3\,e^{11\Lambda q}$ in the same way. By the defintion of $\ell$ in (\ref{eq:subor-l}), we have
\begin{equation}\label{eq:lqlow}
K^6\,e^{22\Lambda q}\ge \|u^{\varphi}\|_{\ell(\varepsilon)}\|v^{\varphi(\varepsilon)}\|_{\ell}=\frac{1}{2\varepsilon}\ge e^{28\Lambda q}
\end{equation}
Therefore, $K>e^{\Lambda q}$ as claim in (\ref{eq:K-lower}). 

To bound $\|u^{\varphi}\|_{\ell}^2$  from below, we need to consider two cases of initial value $\varphi$.\\

\noindent{\bf Case I:} Assume $\varphi$ satisfies
	\begin{equation}\label{eq:phi-case1}
	\|\big(\widetilde A({q})-\frac{{\rm Tr\,}\widetilde A(q)}{2}\cdot I\big)\cdot {\widetilde X}\|\ge e^{-\frac{1}{4}\Lambda q}.
	\end{equation}
	
By (\ref{eq:JZ2.9-1}), for any $e^{\frac{1}{2}\Lambda q}\le k\le K\le N_q,$ we have
	\begin{eqnarray*}
		\|\widetilde A^k({q})\cdot {\widetilde X}\|
		&\ge & \frac{\rho^k-\rho^{-k}}{\rho-\rho^{-1}}\cdot\|\Big(
		\widetilde A(q)-\frac{{\rm Tr\,}\widetilde A(q)}{2}\cdot I\Big){\widetilde X}\|-\frac{\rho^k+\rho^{-k}}{2}\cdot \|{\widetilde X}\| \\
		&\ge& \frac{1}{3}k\cdot e^{-\frac{1}{4}\Lambda q}-C_1 \\
		&\ge& 3,
	\end{eqnarray*}
provided $e^{\frac{1}{4}\Lambda q}>3(C_1+3)$.

By (\ref{eq:A(kq)-Ak(q)}), we have then
	$$\|{\widetilde A}({kq})\cdot {\widetilde X}\|\ge
	\|\widetilde A^k({q})\cdot {\widetilde X}\|-\|\Big({\widetilde A}({kq})-\widetilde A^k({q})\Big)\cdot {\widetilde X}\|
	\ge 2. $$
By (\ref{eq:conjA-A}) and (\ref{eq:rexpbeta}), for $ e^{\frac{1}{4}\Lambda q}\le k\le K$, we have
	\begin{equation}\label{eq:Akq4}
\|{ A}({kq}) { X}\|=|r^{-1}(kq,0)|\cdot \|T_{kq}{\widetilde A}({kq})T_0^{-1} {X}\|=|r^{-1}(kq,0)|\cdot \|{\widetilde A}({kq}){\widetilde X}\|
\ge1
\end{equation}
Therefore,
	$$
	\|u^{\varphi}\|_{\ell}^2
	\ge \frac{1}{2} \sum_{n=1}^{[\ell]-1}\|A_n\cdot { X}\|^2
	\ge \frac{1}{2}\ \ \sum_{e^{\frac{1}{4}\Lambda q}\le k\le K}\|{ A}({kq})\cdot { X}\|^2
	\ge \frac{1}{2}(K-e^{\frac{1}{4}\Lambda q})>e^{\frac{1}{2}\Lambda q}
	$$
	
\noindent{\bf Case II:} Assume $\varphi$ satisfies
	\begin{equation}\label{eq:phi-case2}
	\|\big(\widetilde A({q})-\frac{{\rm Tr\,}\widetilde A(q)}{2} I\big)\cdot {\widetilde X}\|< e^{-\frac{1}{4}\Lambda q},
	\end{equation}
	By (\ref{eq:JZ2.9-1}), for any $1\le k\le e^{\frac{1}{5}\Lambda q}<N_q$ we get
	\begin{eqnarray*}
		\|\widetilde A^k({q})\cdot {\widetilde X}\|
		&\ge & \frac{\rho^k+\rho^{-k}}{2}\cdot \|{\widetilde X}\|-\frac{\rho^k-\rho^{-k}}{\rho-\rho^{-1}}\cdot\|\Big(
		\widetilde A(q)-\frac{{\rm Tr\,}\widetilde A(q)}{2} I\Big){\widetilde X}\| \\
		&\ge& \frac{1}{2}-C_1k\cdot e^{-\frac{1}{4}\Lambda q}\\
		&\ge& \frac{1}{3},
	\end{eqnarray*}
provided $e^{\frac{1}{20}\Lambda q}>6C_1$.

By (\ref{eq:A(kq)-Ak(q)}), we have
	$$\|{\widetilde A}({kq})\cdot {\widetilde X}\|\ge
	\|\widetilde A^k({q})\cdot {\widetilde X}\|-\|({\widetilde A}({kq})-\widetilde A^k({q}))\cdot {\widetilde X}\|
	\ge \frac{1}{4}.$$
	
	By (\ref{eq:conjA-A}) and (\ref{eq:rexpbeta}), for any $1\le k\le e^{\frac{1}{5}\Lambda q}<N_q$, we have
	\begin{equation}\label{eq:Akq5}
	\|{ A}({kq}) { X}\|=|r^{-1}(kq,0)|\cdot \|T_{kq}{\widetilde A}({kq})T_0^{-1} {X}\|=|r^{-1}(kq,0)|\cdot \|{\widetilde A}({kq}){\widetilde X}\|
	\ge \frac{1}{5}
	\end{equation}
	Therefore,
	$$
	\|u^{\varphi}\|_{\ell}^2
	\ge \frac{1}{2} \sum_{n=1}^{[\ell]-1}\|A_n\cdot { X}\|^2
	\ge \frac{1}{2}\ \ \sum_{1\le k\le e^{\frac{1}{5}\Lambda q}} \|{ A}({kq})\cdot { X}\|^2
	\ge \frac{1}{50}e^{\frac{1}{5}\Lambda q}\ge e^{\frac{1}{6}\Lambda q}.
	$$	

\qed

\subsection{The refined estimate on the preimage of ${\cal P}_n(\R)$.}\label{A:JM}
Let ${\cal P}_n(\R)$ denote the polynomials over $\R$ of exact degree
$n$. Let  the class ${\cal P}_{n;n}(\R)$ be elements in ${\cal P}_n(\R)$ with $n$ distinct real zeros. The following proposition was proved in Theorem 6.1,\cite{JMa12Gafa}:
\begin{prop}\label{jm}
Let $p \in {\cal P}_{n;n}(\R)$ with $y_1 < \dots < y_{n-1}$  the
local extrema of $p$. Let
\begin{equation} \label{eq_zeta}
\zeta(p) := \min_{1 \leq j \leq n-1} |p(y_{j})|
\end{equation} and $0 \leq  a < b$.  Then,
\begin{eqnarray}
|p^{-1}(a,b)|  \leq 2 \textrm{diam} (z(p-a)) \max\Big\{ \frac{b-a}{\zeta(p)+a}, \big(\frac{b-a}{\zeta(p)+a}\big)^\frac{1}{2} \Big\}
\end{eqnarray}

where $z(p)$ is the zero set of $p$ and $|\cdot |$ denotes the
Lebesgue measure.
\end{prop}

\section*{Acknowledgement}
The authors would like to thank Svetlana Jitomirskaya for reading  the early manuscript and valuable  suggestions.
The authors would also like to thank Ilya Kachkovskiy for useful comments.
R. H. and F. Y. would like to thank the Institute for Advanced Study, Princeton, for its hospitality during the 2017-18 academic year.
R. H. and F. Y. were supported in part by NSF grant DMS-1638352.
Research of S. Z. was supported in part by NSF grant DMS-1600065.

\bibliographystyle{amsplain}

\vspace{1cm}
Rui Han, 

School of Math, Institute for Advanced Study.

E-mail address:  rhan@ias.edu

\vspace{1cm}
Fan Yang,

School of Math, Institute for Advanced Study.

E-mail address: yangf@ias.edu

\vspace{1cm}
Shiwen Zhang, 

Dept. of Math., Michigan State University. 

E-mail address: zhangshiwen@math.msu.edu

\end{document}